\documentclass[a4paper]{scrartcl}
\usepackage{amsmath,amssymb,url,graphicx, verbatim}
\usepackage{varioref}
\usepackage{relsize}
\newtheorem{theorem}{Theorem}[section]
\newtheorem{lemma}[theorem]{Lemma}

\numberwithin{equation}{section}
\newcommand{\extcite}[2]{\cite[#2]{#1}}

\begin{document}

\begin{center}
 \LARGE \textbf{LIMIT THEOREMS FOR MULTIVARIATE LACUNARY SYSTEMS}
\end{center}

\vspace{3ex}

\begin{center}
 THOMAS L\"OBBE
\footnote{The results are part of the author's PhD thesis supported by IRTG 1132, University of Bielefeld}
\end{center}

\vspace{3ex}

\begin{abstract}
ABSTRACT.  Lacunary function systems of type $(f(M_nx))_{n\geq 1}$ for periodic functions $f$ and sequences of fast-growing matrices $(M_n)_{n\geq 1}$ exhibit many properties of independent random variables like satisfying the Central Limit Theorem or the Law of the Iterated Logarithm. It is well-known that this behaviour depends on number theoretic properties of $(M_n)_{n\geq 1}$ as well as analytic properties of $f$. Classical techniques are essentially based on Fourier analysis making it almost impossible to use a similar approach in the multivariate setting. Recently Aistleitner and Berkes introduced a new method proving the Central Limit Theorem in the one-dimensional case by approximating $\sum_{n}f(M_nx)$ by a sum of piecewise constant periodic functions which form a martingale differences sequence and using a Berry-Esseen type inequality. Later this approach was used to show the Law of the Iterated Logarithm by a consequence of Strassen's almost sure invariance principle. In this paper we develop this method to prove the Central Limit Theorem and the  Law of the Iterated Logarithm in the multidimensional case.
\end{abstract}

\vspace{5ex}

\section{Introduction}

\subsection*{Discrepancy and Uniform Distribution}

A sequence of vectors $(x_n)_{n\geq 1}=(x_{n,1},\ldots,x_{n,d})_{n\geq 1}$ of real numbers in $[0,1)^d$ is called \textit{uniformly distributed modulo one} if
\begin{equation}
 \label{81}
 \lim_{N\to\infty}\frac{1}{N}\sum_{n=1}^N\mathbf{1}_{\mathcal{A}}(x_n)=\lambda(\mathcal{A})
\end{equation}
for any axis-parallel box $\mathcal{A}\subset [0,1)^d$ where $\mathbf{1}_{\mathcal{A}}$ denotes the indicator function on the set $\mathcal{A}$ and $\lambda$ denotes the Lebesgue-measure on $[0,1)^d$.
The discrepancy resp. the star discrepancy of the first $N$ elements of $(x_n)_{n\geq 1}$ is defined by
\begin{equation}
 \label{65}
 \begin{aligned}
 D_N(x_1,\ldots,x_N)\:\:\: & = \:\:\: \sup_{\mathcal{A}\in\mathcal{B}}\left|\frac{1}{N}\sum_{n=1}^N\mathbf{1}_{\mathcal{A}}(x_n)-\lambda(\mathcal{A})\right|,\\
 D^*_N(x_1,\ldots,x_N)\:\:\: & =\:\:\: \sup_{\mathcal{A}\in\mathcal{B}^*}\left|\frac{1}{N}\sum_{n=1}^N\mathbf{1}_{\mathcal{A}}(x_n)-\lambda(\mathcal{A})\right|,
 \end{aligned}
\end{equation}
where $\mathcal{B}$ denotes the set of all axis-parallel boxes $\mathcal{A}=\prod_{i=1}^d[\alpha_{i},\beta_{i})\subset [0,1)^d$ and furthermore
$\mathcal{B}^*$ denotes the set of all axis-parallel boxes $\mathcal{A}=\prod_{i=1}^d[0,\beta_{i})\subset [0,1)^d$ with one corner in $0$. 
It is well-known that (\ref{81}) is equivalent to $D_N(x_1,\ldots,x_N)\to 0$ resp. $D^*_N(x_1,\ldots,x_N)\to 0$ for $N\to\infty$. By a classical result of Weyl \cite{W16} it is known that for any increasing sequence $(M_n)_{n\geq 1}$ of positive integers the sequence $(\langle M_nx\rangle)_{n\geq 1}$, where $\langle\cdot\rangle$ denotes the fractional part, is uniformly distributed modulo one for almost all $x\in [0,1)$. This result naturally extends to the multidimensional case. Sequences with vanishing star-discrepancy have applications in the theory of numerical integration. The connection is established by the Koksma-Hlawka inequality (see \cite{DT97}) which states that for any sequence of vectors $(x_n)_{n\geq 1}\subset [0,1)^d$ we have
\begin{equation}
 \label{217}
 \left|\frac{1}{N}\sum_{n=1}^Nf(x_n)-\int_{[0,1)^d}f(x)\,dx\right|\leq D_N^*(x_1,\ldots,x_N)\cdot V_{HK}(f)
\end{equation}
for any function $f$ on $[0,1)^d$ where $V_{HK}$ denotes the total variation in the sense of Hardy and Krause. Thus the integral can be approximated by the mean of the values which some points have under $f$ where the approximation error is given by the total variation of $f$ and the star-discrepancy of the points. Although (\ref{217}) is interesting form a theoretical point of view, it is of little use in practice. In general the total variation is more difficult to compute than the integral. But nevertheless, it becomes evident that sequences of points with low discrepancy give small approximation errors. Therefore we are not only interested in sequences such that the star-discrepancy tends to 0, but also in the speed of convergence.

\subsection*{Lacunary sequences}
Let $(M_n)_{n\geq 1}$ be a sequence of non-singular integer-valued $d\times d$-matrices satisfying a Hadamard gap type condition of the form
\begin{equation}
\label{17}
||M_{n+k}^Tj||_{\infty}\geq q^k||M_n^T||_{\infty}
\end{equation}
for all $j\in\mathbb{Z}^d\backslash\{0\}$, $n\in\mathbb{N}$, $k\geq \log_q(||j||_{\infty})$ and some absolute constant $q>1$. Here $A^T$ denotes the transpose of a matrix $A$.  Since this extends the definition of lacunary sequences for $d=1$ to the multivariate case we call this system a multivariate lacunary sequence satisfying a Hadamard gap condition.
For $d=1$ and some sequence $(M_n)_{n\geq 1}$ satisfying (\ref{17}) Salem and Zygmund \cite{SZ47} proved that for any sequence of integers $(a_n)_{n\geq 1}$ with
\begin{equation*}
 a_N=o(A_N) \quad \textnormal{for} \quad A_N=\frac{1}{2}\left(\sum_{n=1}^Na_n^2\right)^{1/2}
\end{equation*}
we have
\begin{equation}
\lim_{N\to\infty}\mathbb{P}\left(\frac{1}{A_N}\sum_{n=1}^Na_n\cos(2\pi M_nx)\leq t \right)=\Phi(t)
\end{equation}
where $\mathbb{P}$ denotes the probability measure induced by the Lebesgue measure on $[0,1)^d$ and $\Phi$ denotes the standard normal distribution, i.e. for all $t\in\mathbb{R}$ we have
\begin{equation*}
 \Phi(t)=\frac{1}{\sqrt{2\pi}}\int_{-\infty}^{t}e^{-\frac{1}{2}y^2}\,dy.
\end{equation*}

Furthermore Weiss \cite{W59} (see also Salem and Zygmund \cite{SZ50}, Erd\H{o}s and G\'al \cite{EG55}) showed that
\begin{equation}
\limsup_{N\to\infty}\frac{\sum_{n=1}^Na_n\cos(2\pi M_nx)}{\sqrt{2A_N^2\log\log(N)}}=1 \quad a.e.
\end{equation}
under the condition
\begin{equation*}
 a_N=o\left(\frac{A_N}{\sqrt{\log(\log(A_N))}}\right).
\end{equation*}
Therefore for lacunary $(M_n)_{n\geq 1}$ the sequence $(a_n\cos(2\pi M_nx))_{n\geq 1}$ shows a behaviour typical for independent, identically distributed random variables. One could ask whether this holds for other periodic
functions as well. The answer is negative in general. By a result of Erd\H{o}s and Fortet (see \cite{K49}) for $f(x)=\cos(2\pi x)+\cos(4\pi x)$ and $M_n=2^n-1$ we have
\begin{equation}
\lim_{N\to\infty}\mathbb{P}\left(\frac{1}{\sqrt{N}}\sum_{n=1}^Nf(M_nx)\leq t\right)=\frac{1}{\sqrt{\pi}}\int_{0}^1\int_{-\infty}^{t|\cos(\pi s)|/2}e^{-u^2}\,duds
\end{equation}
and
\begin{equation}
\limsup_{N\to\infty}\frac{\sum_{n=1}^Nf(M_nx)}{\sqrt{N\log(\log(N))}}=2\cos(\pi x) \quad a.e.
\end{equation}
Thus neither the Central Limit Theorem nor the Law of the Iterated Logarithm is satisfied. This result was later generalized by Conze and Le Borgne \cite{CB11} (see also \cite{A13a} for further information).
On the other hand Kac \cite{K46} showed that any one-periodic function $f:\mathbb{R}\to\mathbb{R}$ of mean zero which is of bounded variation on $[0,1)$ or Lipschitz-continuous satisfies
\begin{equation}
\lim_{N\to\infty}\mathbb{P}\left(\frac{1}{\sqrt{N}}\sum_{n=1}^Nf(2^nx)\leq t\sigma\right)=\Phi(t)
\end{equation}
if
\begin{equation}
\sigma^2=\mathbb{E}[f]+2\sum_{n=1}^{\infty}\mathbb{E}[f(x)f(2^nx)]\neq 0.
\end{equation}
Furthermore Maruyama \cite{M50} and Izumi \cite{I51} proved
\begin{equation}
 \limsup_{N\to\infty}\frac{\sum_{n=1}^Nf(2^nx)}{\sqrt{2N\log(\log(N))}}=\sigma \quad \textnormal{a.e.}
\end{equation}
This illustrates that the behaviour of $(f(M_nx))_{n\geq 1}$ does not only depend on the speed of growth of $(M_n)_{n\geq 1}$ but also on number theoretic properties of the sequence $(M_n)_{n\geq 1}$.
Later on the Central Limit Theorem was shown for more general lacunary sequences. By a result of Gaposhkin \cite{G66}
\begin{equation}
 \label{1}
 \lim_{N\to\infty}\mathbb{P}\left(\sum_{n=1}^Nf(M_nx)\leq t\sigma_N\right)=\Phi(t)
\end{equation}
holds for sequences $(M_n)_{n\geq 1}$ satisfying
\begin{equation}
 \label{2}
 \sigma_N^2=\mathlarger{\mathlarger{\int}}_0^1\left(\sum_{i=1}^Nf(M_nx)\right)^2\,dx\geq CN,
\end{equation}
for an absolute constant $C>0$ and one of the following conditions
\begin{itemize}
 \item $\frac{M_{n+1}}{M_n}\in\mathbb{N}$, \quad \textnormal{for all } $n\in\mathbb{N}$,
 \item $\lim_{n\to\infty}\frac{M_{n+1}}{M_n}=\theta$, \quad such that $\theta^r$ irrational for all $r\in\mathbb{N}$.
\end{itemize}
Takahashi \cite{T61} showed (\ref{1}) for $M_{n+1}/M_{n}\to\infty$ and $\alpha$-Lipschitz-continuous functions.
The connection between the Central Limit Theorem and the number of solutions of certain Diophantine equations is due to Gaposhkin \cite{G70}.
Consider the linear Diophantine equation
\begin{equation*}
 aj\pm a'j'=\nu
\end{equation*}
for fixed integers $j,j',\nu$. In general the set of solutions consists of all pairs of integers $a,a'$ such that equality holds but we restrict ourselves to those solutions with $a=M_n$ and $a'=M_{n'}$ for $n,n'\in\mathbb{N}$ and rather regard the indices $n,n'$ as solutions of this equation.
The Central Limit Theorem holds for lacunary sequences $(M_n)_{n\geq 1}$ satisfying (\ref{2}), if for any fixed $j,j',\nu$ the number of solutions of the Diophantine equation
\begin{equation}
\label{61}
 M_nj\pm M_{n'}j'=\nu
\end{equation}
is bounded by an absolute constant $C_{j,j'}>0$ which is independent of $\nu$. Observe that ``nice'' periodic functions can be approximated by trigonometric polynomials very well. Thus because of the product-to-sum identities of trigonometric functions the behaviour of the moments of $\sum f(M_nx)$ depends on the number of solutions of Diophantine equations of certain length.
Recently Aistleitner and Berkes \cite{AB10} improved this result: For a lacunary sequence $(M_n)_{n\geq 1}$ satisfying the Hadamard gap condition set
\begin{equation}
 \label{22}
 \begin{aligned}
  L(N,G,\nu) \:\:\:  = & \:\:\: |\{1\leq n,n'\leq N:\\
  &\:\:\: \exists j,j'\in\mathbb{Z}^d,1\leq ||j||_{\infty}|,||j'||_{\infty}\leq G,M_n^Tj\pm M_{n'}^Tj'=\nu\}|,\\
  L^*(N,G,\nu) \:\:\: = & \:\:\: |\{1\leq n,n'\leq N,n\neq n':\\
  &\:\:\: \exists j,j'\in\mathbb{Z}^d,1\leq ||j||_{\infty}|,||j'||_{\infty}\leq G,M_n^Tj\pm M_{n'}^Tj'=\nu\}|,\\
  L(N,G) \:\:\: = & \:\:\: \sup_{\nu\neq 0}L(N,G,\nu).
 \end{aligned}
\end{equation}
and for all $N\geq 1$, $G\geq 1$ and $\nu\in\mathbb{Z}$. Let $f:\mathbb{R}\to\mathbb{R}$ be some function of finite total variation which is one-periodic and satisfies $\mathbb{E}[f]=0$ as well as (\ref{2}) for some lacunary sequence satisfying the Hadamard gap condition (\ref{17}).
Aistleitner and Berkes showed that if for any fixed $G\geq 1$ we have $L(N,G)=o(N)$ for $N\to\infty$ then (\ref{1}) holds.

\subsection*{Law of the Iterated Logarithm for the discrepancy of lacunary point sets}

The Law of the Iterated Logarithm for the discrepancy of an one-dimensional lacunary point set was shown by Philipp \cite{P75}. He proved
\begin{equation*}
 \frac{1}{4\sqrt{2}}\leq \limsup_{N\to\infty}\frac{ND_N(M_1x,\ldots,M_Nx)}{\sqrt{2N\log(\log(N))}}\leq C \quad \textnormal{a.e.}
\end{equation*}
 where the constant $C>0$ depends on $q$ only. This corresponds to the Chung-Smirnov Law of the Iterated Logarithm, that is
\begin{equation}
 \label{220}
 \limsup_{N\to\infty}\frac{ND_N(\xi_1,\ldots,\xi_N)}{\sqrt{2N\log(\log(N))}}=\frac{1}{2} \quad \textnormal{a.s.}
\end{equation}
for any sequence of independent, identically distributed non-degenerate random variables $(\xi_n)_{n\geq 1}$ in $[0,1)$ with $\mathbb{E}[\xi_1]=0$ and $\mathbb{E}[\xi_1^2]=1$. For sequences of type $(M_n)_{n\geq 1}=(\theta^n)_{n\geq 1}$ for $\theta>0$ the precise value of the Law of the Iterated Logarithm was determined by Fukuyama \cite{F08}, i.e. for a.e. $x$ we have
\begin{equation*}
 \limsup_{N\to\infty}\frac{ND_N(\theta^1x,\ldots, \theta^Nx)}{\sqrt{2N\log(\log(N))}}=
\begin{cases}
 \sqrt{42}/9, & \textnormal{if }\theta=2,\\
 \frac{\sqrt{(\theta+1)\theta(\theta-2)}}{2\sqrt{(\theta-1)^3}}, & \textnormal{if }\theta\geq 4\textnormal{ is an even integer},\\
 \frac{\sqrt{(\theta+1)}}{2\sqrt{\theta-1}}, & \textnormal{if }\theta\geq 3\textnormal{ is an odd integer},\\
 1/2, & \textnormal{if } \theta^r\not\in\mathbb{Q} \textnormal{ for all }r\in\mathbb{N}.
\end{cases}
\end{equation*}
Therefore the probabilistic analogy is not complete. The precise value depends sensitively on number theoretic properties of the sequence $(M_n)_{n\geq 1}$, mainly on the number of non-trivial solutions of the Diophantine equations $M^T_nj\pm M^T_{n'}j'=0$.\\
Aistleitner \cite{A10} used the method applied in \cite{AB10} to prove the Law of the Iterated Logarithm for function systems $(f(M_nx))_{n\geq 1}$ as well as for the discrepancy $D_N(M_1x,\ldots,M_Nx)$ for point sets defined by a lacunary sequence $(M_n)_{n\geq 1}$ satisfying the Hadamard gap condition (\ref{17}) if for any fixed $G\geq 1$ we have $\max(L(N,G),L^*(N,G,0))=\mathcal{O}(N/(\log(N))^{1+\varepsilon})$ for $N\to\infty$. Later Aistleitner, Fukuyama and Furuya \cite{AFF13} improved this result by proving sufficiency of $L(N,G)=\mathcal{O}(N/(\log(N))^{1+\varepsilon})$ for the Law of the Iterated Logarithm for lacunary function systems and $L^*(N,G,0)=o(N)$ in addition to the former condition for the Law of the Iterated Logarithm for discrepancy of lacunary point sets.\\

\subsection*{Functions in several variables}

The Central Limit Theorem for lacunary sequences $(M_n)_{n\geq 1}$ of $d\times d$-matrices satisfying (\ref{17}) was proved by Conze, Le Borgne and Roger \cite{CBR12}. There it was shown that the Central Limit Theorem holds if the sequence is satisfying a strong number theoretic condition, i.e. there is an absolute constant $C$ such that for any integers $G$ and $N$ the following condition holds: For $2s$ integers $1\leq n_1\leq n_1'<n_2\leq n_2'<\cdots<n_s\leq n_s'\leq N$ with $n_{k+1}\geq n_k'+C\log_q(G)$ for $k\in\{1,\ldots,s-1\}$ and vectors $j_1,j_1',\ldots,j_s,j_s'$ with $||j_k||_{\infty},||j_k'||_{\infty}\leq G$ for $k\in\{1,\ldots, s\}$ we have
\begin{equation*}
 M^T_{n_s}j_s+M^T_{n_s'}j_s'\neq 0 \quad \Longrightarrow \quad \sum_{k=1}^sM^T_{n_k}j_k+M^T_{n_k'}j_k'\neq 0.
\end{equation*}
Such condition for example is satisfied in the product case, i.e. there exists a sequence of matrices $(A_n)_{n\geq 1}$ with $M^T_n=A^T_1\cdots A^T_n$ for all $n\in\mathbb{N}$.

\subsection*{Main results}

The methods used in early results in this area are based on substantial use of Fourier analysis such as bounding the size of Fourier coefficients, the tails of Fourier series etc.
In \cite{AB10} a new method effectively reducing the use of Fourier analysis was introduced which was used in \cite{A10} resp. \cite{AFF13} to show the Law of the Iterated Logarithm.

We adopt this method to prove the Central Limit Theorem and the Law of the Iterated Logarithm for lacunary sequences $(M_n)_{n\geq 1}$ satisfying a Hadamard gap condition (\ref{17}) under some weak conditions on the number of solutions of the Diophantine equation $M_n^Tj\pm M_{n'}^Tj'=\nu$.
The Central Limit Theorem for multivariate lacunary systems reads as follows
\begin{theorem}[Central Limit Theorem]
\label{3}
 Let $(M_n)_{n\geq 1}$ be a lacunary sequence of non-singular $d\times d$-matrices satisfying the Hadamard gap condition (\ref{17}).
 Furthermore assume that $L(G,N)=o(N)$ for any fixed $G\geq 2$.
Let $f\in L^2(\mathbb{R}^d,\mathbb{R})$ be a bounded, periodic function with mean zeros which is of finite total variation in the sense of Hardy and Krause. Assume that there exists some absolute constant $C>0$ such that
\begin{equation}
\label{21}
 \sigma_N^2:=\mathlarger{\mathlarger{\int}}_{[0,1)^d}\left(\sum_{n=1}^Nf(M_nx)\right)^2\,dx\geq C\cdot N.
\end{equation}
for any $N\geq 1$. Then for all $t\in\mathbb{R}$ we have
\begin{equation}
 \label{138}
 \lim_{N\to\infty}\left|\mathbb{P}\left(\sum_{n=1}^Nf(M_nx)\leq t\sigma_N\right)-\Phi(t)\right|=0.
\end{equation}
If furthermore for some $0<\beta<1$ we have $L(N,G_N)=\mathcal{O}(N^{\beta})$ for any sequence $(G_N)_{N\geq 1}$ with $G_N\leq dN$ for all $N\geq 1$
then for all $t\in\mathbb{R}$ and sufficiently large $N$ we get
\begin{equation}
 \label{139}
 \left|\mathbb{P}\left(\sum_{n=1}^Nf(M_nx)\leq t\sigma_N\right)-\Phi(t)\right|\leq C\frac{d^{1/5}\log(N)^{3/5}+\log(d)\log(N)}{N^{\min(1/8,(1-\beta)/5)}}
\end{equation}
with some absolute constant $C>0$ which only depends only on $q$.
\end{theorem}
Under a slightly stronger condition on the number of solutions of the Diophantine equation we also obtain

\begin{theorem}[Law of the Iterated Logarithm]
\label{16}
Let $(M_n)_{n\geq 1}$ be a sequence of non-singular $d\times d$-matrices satisfying the Hadamard gap condition (\ref{17}).
Furthermore assume that for any fixed $G\geq 1$ and some $\varepsilon>0$ we have $L(N,G)=\mathcal{O}(N/(\log N)^{1+\varepsilon})$.
Let $f\in L^2(\mathbb{R}^d,\mathbb{R})$ be a bounded, periodic function with mean zeros which is of finite total variation in the sense of Hardy and Krause.
Additionally, let $f$ and $(M_n)_{n\geq 1}$ be given such that for
\begin{equation*}
\sigma_N^2:=\mathlarger{\mathlarger{\int}}_{[0,1)^d}\left(\sum_{n=1}^Nf(M_nx)\right)^2\,dx
\end{equation*}
there exists $\Sigma_{f,M_n}>0$ with
\begin{equation}
 \lim_{N\to\infty}\frac{\sigma_N^2}{N}=\Sigma_{f,M_n}.
\end{equation}
Then we have
\begin{equation}
 \limsup_{N\to\infty}\frac{|\sum_{k=1}^Nf(M_nx)|}{\sqrt{2N\log(\log(N))}}=\sqrt{\Sigma_{f,M_n}} \quad \textnormal{a.e.}
\end{equation}
\end{theorem}

We now state a version of the Law of the Iterated Logarithm for the discrepancy of point sets defined by multivariate lacunary sequences with not too many non-trivial solutions of the Diophantine equation:

\begin{theorem}
\label{63}
Let $(M_n)_{n\geq 1}$ be a lacunary sequence of non-singular $d\times d$-matrices satisfying the Hadamard gap condition (\ref{17}).
Assume that $L(N,G)=\mathcal{O}(N/(\log{N})^{1+\varepsilon})$ and furthermore $L^*(N,G,0)=o(N)$.
Then the discrepancy of $(M_nx)_{n\geq 1}$ resp. the star discrepancy satisfies the Law of the Iterated Logarithm, i.e.
\begin{equation}
 \label{66}
\limsup_{N\to\infty}\frac{ND_N(M_1x,\ldots,M_Nx)}{\sqrt{2N\log(\log(N))}}=\limsup_{N\to\infty}\frac{ND^*_N(M_1x,\ldots,M_Nx)}{\sqrt{2N\log(\log(N))}}=\frac{1}{2} \quad \textnormal{a.e.}
\end{equation}
\end{theorem}

The main idea in the proof of the Central Limit Theorem is to apply a Theorem due to Heyde and Brown \cite{HB70} which ensures the Central Limit Theorem for martingale differences sequences satisfying certain moment conditions resp. a consequence of Strassen's almost sure invariance principle \cite{S67} for which we get the Law of the Iterated Logarithm under similar moments conditions.
The elements of the martingale differences are defined by sums of the form $\sum_{n\in\Delta_k}\varphi_n(x)$ where $\varphi_n(x)$ is a piecewise constant function approximating $f(M_nx)$. The blocks $\Delta_k$ form a sequence of growing blocks which decomposes the set of natural numbers except for small gaps between consecutive blocks.
The filtration is defined by a sequence of $\sigma$-fields which are generated by a decomposition of $[0,1)^d$ into ``dyadic'' blocks, i.e. their side lengths are negative powers of $2$ where the exponents depend on the magnitude of the ``frequencies'' $M_n$ in the corresponding block in a certain manner.
The Central Limit Theorem is ensured by a Berry-Esseen type inequality which gives an upper bound on the mutual distance between the distribution function of the normalized sums of the martingales and the distribution function of a standard normal distributed random variable which only depends on second and fourth moments conditions.
In fact only for an upper bound of the conditional variances a condition on the number of solutions of the Diophantine equations are necessary since all other moments for which we need upper bounds can be estimated in a different way.


\section{Preliminaries}

In this section we repeat some basic results on periodic functions of finite total variation resp. lacunary sequences which are going to be used in the subsequent sections.

For some integer $d\geq 1$ set $I=\{1,\ldots,d\}$.
We now introduce the total variation in the sense of Hardy and Krause for periodic functions on $\mathbb{R}^d$.
Let $f:\mathbb{R}^d\to\mathbb{R}$ be some periodic function, i.e. $f$ satisfies $f(x+z)=f(x)$ for all $x\in\mathbb{R}^d$ and $z\in\mathbb{Z}^d$. For some subset $J\subseteq I$ and points $a,b\in [0,1)^{|J|}$ with $a_i\leq b_i$ for all $i\in J$ and some $z\in [0,1)^{|I\backslash J|}$ define
\begin{equation*}
 \Delta_J(f,a,b,z)=\sum_{\delta\in \{0,1\}^{|J|}}(-1)^{\sum_{i\in J}\delta_i}f(c_{\delta})
\end{equation*}
where $c_{\delta}=(c_{\delta,1},\ldots,c_{\delta,d})$ is defined by $c_{\delta,i}=\delta_i a_i+(1-\delta_i)b_i$ for $i\in J$ and $c_i=z_i$ for $i\notin J$. A finite set $\mathcal{Y}_i=\{y_1\ldots,y_{m(i)}\}\subset [0,1)$ with $0=y_1<\ldots<y_{m(i)}<1$ for some positive integer $m(i)$ is called a ladder. A multidimensional ladder on $[0,1)^d$ has the form $\mathcal{Y}=\prod_{i\in I}\mathcal{Y}_i$. For a multidimensional ladder $\mathcal{Y}$, a subset $J\subseteq I$ and $z\in[0,1)^d$ set $\mathcal{Y}_{J,z}=\prod_{i\in J}\mathcal{Y}_i\times \prod_{i\notin J}\{z_i\}$.
For $y\in\mathcal{Y}_{J,z}$ define $y_+\in [0,1)^{|J|}\times \prod_{i\notin J}\{z_i\}$ such that $y_{+,i}$ is the successor of $y_i$ in $\mathcal{Y}_i$ resp. $1$ if $y_i$ is the largest element in $\mathcal{Y}_i$. Then we define the variation of $f$ over $\mathcal{Y}_J$ by
\begin{equation*}
 V_{\mathcal{Y}_{J,z}}(f)=\sum_{y\in\mathcal{Y}_{J,z}}|\Delta_J(f,y,y_+,z)|.
\end{equation*}
Denote the set of all ladders $\mathcal{Y}_{J,z}$ by $\mathbb{Y}_{J,z}$. Then the total variation of $f$ over $[0,1)^{|J|}$ is defined by
\begin{equation*}
 V_{J}(f)=\sup_{z\in[0,1)^{|I|}}\sup_{\mathcal{Y}_{J,z}\in\mathbb{Y}_{J,z}}V_{\mathcal{Y}_{J,z}}(f).
\end{equation*}
The total variation of $f$ on $[0,1)^d$ in the sense of Hardy and Krause is
\begin{equation*}
 V_{HK}(f)=\sum_{J\subseteq I,J\neq \emptyset}V_{J}(f).
\end{equation*}
A function $f$ is called to be of finite total variation if $V_{HK}(f)<\infty$.

The following Lemma was proved in \cite{Z68}:

\begin{lemma}
\label{40}
 Let $f(x)=\sum_{j\in\mathbb{Z}^d\backslash\{0\}}a_j\cos(2\pi\langle j,x\rangle)+b_j\sin(2\pi\langle j,x\rangle)$ be a periodic function of finite total variation in the sense of Hardy and Krause. Then we have
\begin{equation*}
|a_j|,|b_j|\leq C\left(\prod_{i\in I,j_i\neq 0}\frac{1}{2\pi |j_i|}\right)V_{\{i\in I: j_i\neq 0\}}(f).
\end{equation*}
for some absolute constant $C>0$ and all $j\in\mathbb{Z}^d\backslash\{0\}$.
\end{lemma}

Observe that hereafter we always write $C$ for some absolute constant which may vary from line to line. Furthermore we always assume that $f\in L^2(\mathbb{R}^d,\mathbb{R})$ is a periodic function of finite total variation in the sense of Hardy and Krause.

For $\Gamma\in\mathbb{N}_0^d$ we denote the $\Gamma$th Dirichlet kernel by
\begin{equation}
\label{301}
D_{\Gamma}(x)=\sum_{j\in\mathbb{Z}^d,|j_i|\leq \Gamma_i}\cos(2\pi\langle j,x\rangle).
\end{equation}
Then the $\Gamma$th partial sum of $f$ is defined by
\begin{equation}
 \label{127}
 \psi_{\Gamma}(x)=\int_{[0,1)^d}f(x+t)D_{\Gamma}(t)\,dt=\sum_{j\in\mathbb{Z}^d,|j_i|\leq\Gamma_i}a_j\cos(2\pi\langle j,x\rangle)+b_j\sin(2\pi\langle j,x\rangle)
\end{equation}
for suitable numbers $a_j,b_j\in\mathbb{R}$ for all $j\in\mathbb{Z}^d$ with $|j_i|\leq \Gamma$ for all $i\in\{1,\ldots,d\}$. Set $\rho_{\Gamma}(x)=f(x)-\psi(x)$. If $\Gamma_i=G$ for all $i\in I$ we simply write $\psi_G(x)$ resp. $\rho_G(x)$. The $G$th Fej\'er mean of $f$ is defined by
\begin{equation}
 \label{92}
 \begin{aligned}
  p_G(x) \:\:\: & = \:\:\: \frac{1}{(G+1)^d}\sum_{\Gamma\in\mathbb{N}_0^d,||\Gamma||_{\infty}\leq G}\psi_{\Gamma}(x)\\
  & = \:\:\: \sum_{j\in\mathbb{Z}^d,|j_i|\leq G}a'_j\cos(2\pi\langle j,x\rangle)+b'_j\sin(2\pi\langle j,x\rangle)
 \end{aligned}
\end{equation}
where
\begin{equation*}
 a'_j=a_j\prod_{i\in I}\frac{G+1-|j_i|}{G+1}, \quad  b'_j=b_j\prod_{i\in I}\frac{G+1-|j_i|}{G+1}
\end{equation*}
for all $j\in\mathbb{Z}^d$ with $||j||_{\infty}\leq G$. Observe that
\begin{equation*}
 p_G(x)=\int_{[0,1)^d}f(x+t)K_G(t)\,dt
\end{equation*}
where $K_G(t)=K_{\Gamma}(t)$ with $\Gamma_i=G$ for all $i\in I$ and $K_{\Gamma}(t)=\prod_{i\in I}K_{\Gamma_i}(t_i)=\prod_{i\in I}K_{G}(t_i)$ is the $d$-dimensional $G$th Fej\'er kernel and $K_G(t_i)$ is the one-dimensional $G$th Fej\'er kernel defined by
\begin{eqnarray*}
 K_G(t_i) & = & \frac{1}{G+1}\sum_{l=0}^{G}\sum_{j_i=-l}^{l}\cos(2\pi\langle j_i,x\rangle)\\
 & = & \frac{1}{G+1}\sum_{l=0}^{G}\frac{\sin(2\pi\langle l+1/2,t_i \rangle)}{\sin(2\pi\langle 1/2,t_i \rangle)}\\
 & = & \frac{1}{G+1}\frac{(\sin(2\pi\langle (G+1)/2,t_i \rangle))^2}{2(\sin(2\pi\langle 1/2,t_i \rangle))^2}\\
 & \geq & 0.
\end{eqnarray*}
Therefore we have
\begin{eqnarray*}
 |p_G(x)| & = & \left|\int_{[0,1)^d}f(x+t)K_G(t)\,dt\right|\\
 & \leq & ||f||_{\infty}\left|\int_{[0,1)^d}K_G(t)\,dt\right|\\
 & \leq & ||f||_{\infty}.
\end{eqnarray*}
Now we define $r_G(x)=f(x)-p_G(x)$.

\begin{lemma}
 \label{94}
 Let $f(x)=\sum_{j\in\mathbb{Z}^d\backslash\{0\}}a_j\cos(2\pi\langle j,x\rangle)+b_j\sin(2\pi\langle j,x\rangle)$ be some periodic function satisfying $V_{HK}(f)\leq 1$. Then there exists some absolute constant $C>0$ such that for any $G\geq d$ the function $r_G$ satisfies
 \begin{equation*}
  ||r_G||_2^2\leq CdG^{-1}.
 \end{equation*}
\end{lemma}

\textit{Proof.} 
Observe that for any $\varepsilon=\varepsilon(d,G)>0$ there exists some trigonometric polynomial $f'$ such that $||f-f'||_2\leq \varepsilon$. Let $p_G'$ be the $G$th Fej\'er mean of $f'$. We obtain $||p_G-p_G'||_2\leq \varepsilon$. Therefore we have $||r_G||_2\leq ||f'-p_g'||_2+C\sqrt{dG^{-1}}$. Thus it is enough to prove the statement of the Lemma for trigonometric polynomials $f$.
Set
\begin{equation}
  \label{83}
 r_G(x)=\sum_{j\in\mathbb{Z}^d\backslash\{0\}}\tilde{a}_j\cos(2\pi\langle j,x\rangle)+\tilde{b}_j\sin(2\pi\langle j,x\rangle)
\end{equation}
where
\begin{equation*}
 \tilde{a}_j=
 \begin{cases}
 a_j, & ||j||_{\infty}>G,\\
 a_j\left(1-\prod_{i\in I}\left(1-\frac{|j_i|}{G+1}\right)\right), & ||j||_{\infty}\leq G
 \end{cases}
\end{equation*}
and $\tilde{b}_j$ is defined analogously.
We have
\begin{equation}
  \label{91}
 ||r_G||_2^2=\sum_{j\in\mathbb{Z}^d\backslash\{0\}}\tilde{a}_j^2+\tilde{b}_j^2=\underbrace{\sum_{0<||j||_{\infty}\leq G}\tilde{a}_j^2+\tilde{b}_j^2}_{(*)}+\underbrace{\sum_{||j||_{\infty}>G}\tilde{a}_j^2+\tilde{b}_j^2}_{(**)}.
\end{equation}
For some given nonempty $J\subseteq I$ set 
\begin{eqnarray*}
 D(G,J) & = & \{j\in\mathbb{Z}^d:1\leq |j_i|\leq G\textnormal{ for } i\in J,j_i=0\textnormal{ for } i\notin J\},\\
 D'(G,J) & = & \{j\in\mathbb{Z}^d:j_i\neq 0 \textnormal{ for }i\in J,j_i=0\textnormal{ for } i\notin J,j\notin D(G,J)\}.
\end{eqnarray*}
To estimate $(*)$ we first by Lemma \ref{40} observe 
\begin{equation*}
 \sum_{0<||j||_{\infty}\leq G}\tilde{a}_j^2+\tilde{b}_j^2 \leq  2\sum_{J\subseteq I,J\neq \emptyset}\sum_{j\in D(G,J)}\left(\prod_{i\in J}\frac{1}{2\pi|j_i|}\right)^2\left(1-\prod_{i\in J}\left(1-\frac{|j_i|}{G+1}\right)\right)^2V_J(f)^2.
\end{equation*}
By definition of $V_{HK}(f)$ it is enough to show
\begin{equation}
 \label{89}
 V(G,J)=2\sum_{j\in D(G,J)}\left(\prod_{i\in J}\frac{1}{2\pi|j_i|}\right)^2\left(1-\prod_{i\in J}\left(1-\frac{|j_i|}{G+1}\right)\right)^2\leq CG^{-1}
\end{equation}
for some absolute constant $C>0$. By decomposing we have
\begin{equation*}
 \begin{aligned}
 V(G,J)=2\sum_{K,K'\subseteq J,K,K'\neq 0}\sum_{j\in D(G,J)} & \prod_{i\in K}\frac{1}{2\pi|j_i|}\frac{|j_i|}{G+1}\prod_{i\in J\backslash K}\frac{1}{2\pi|j_i|}\left(1-\frac{|j_i|}{G+1}\right)\\
 \cdot & \prod_{i\in K'}\frac{1}{2\pi|j_i|}\frac{|j_i|}{G+1}\prod_{i\in J\backslash K'}\frac{1}{2\pi|j_i|}\left(1-\frac{|j_i|}{G+1}\right).
 \end{aligned}
\end{equation*}
Thus we get
\begin{equation*}
 V(G,J)=2\sum_{K,K'\subseteq J,K,K'\neq 0}W_1(K,K')\cdot W_2(K,K')\cdot W_3(K,K')\cdot W_4(K,K')
\end{equation*}
where
\begin{eqnarray*}
 W_1(K,K') & = & \prod_{i\in K\cap K'}\frac{1}{(2\pi)^2}\sum_{j_i=-G}^G\frac{1}{(G+1)^2}\leq \prod_{i\in K\cap K'}\frac{1}{4G},\\
 W_2(K,K') & = & \prod_{i\in K\cap J\backslash K'}\frac{1}{(2\pi)^2}\sum_{j_i=-G}^G\frac{1}{G+1}\left(\frac{1}{|j_i|}-\frac{1}{G+1}\right)\leq \prod_{i\in K\cap J\backslash K'}\frac{\log(G)}{4G},\\
 W_3(K,K') & = & \prod_{i\in J\backslash K\cap K'}\frac{1}{(2\pi)^2}\sum_{j_i=-G}^G\left(\frac{1}{|j_i|}-\frac{1}{G+1}\right)\frac{1}{G+1}\leq \prod_{i\in J\backslash K\cap K'}\frac{\log(G)}{4G},\\
 W_4(K,K') & = & \prod_{i\in J\backslash K\cap J\backslash K'}\frac{1}{(2\pi)^2}\sum_{j_i=-G}^G\left(\frac{1}{|j_i|}-\frac{1}{G+1}\right)^2\leq \prod_{i\in J\backslash K\cap J\backslash K'}\frac{1}{4}.
\end{eqnarray*}
Since $K\cap K'\neq\emptyset$ or $K\cap J\backslash K'\neq\emptyset$ and $J\backslash K\cap K'\neq\emptyset$ we conclude
\begin{equation}
 \label{90}
 V(G,J)\leq \sum_{K,K'\subseteq J,K,K'\neq 0}\frac{C}{4^{|J|}G}\leq CG^{-1}
\end{equation}
for some absolute constant $C>0$ and therefore (\ref{89}) is verified.
We now estimate $(**)$. By Lemma \ref{40} we have
\begin{eqnarray*}
 \sum_{||j||_{\infty}>G}\tilde{a}_j^2+\tilde{b}_j^2 & \leq & \sum_{J\subseteq I,J\neq\emptyset}\sum_{j\in D'(G,J)}a_j^2+b_j^2\\
 & \leq & 2\sum_{J\subseteq I,J\neq\emptyset}\sum_{j\in D'(G,J)}\left(\prod_{i\in J}\frac{1}{(2\pi|j_i|)^2}\right)V_J(f)^2.
\end{eqnarray*}
We furthermore for some nonempty $J\subseteq I$ get
\begin{eqnarray*}
 \sum_{j\in D'(G,J)}\left(\prod_{i\in J}\frac{1}{(2\pi|j_i|)^2}\right) & \leq &  \sum_{l=1}^{|J|}{|J|\choose l}\left(2\sum_{j=1}^G\frac{1}{(2\pi j)^2}\right)^{|J|-l}\cdot \left(2\sum_{j=G+1}^{\infty}\frac{1}{(2\pi j)^2}\right)^{l}\\
 & \leq & \sum_{l=1}^{|J|}{|J|\choose l}\left(\frac{1}{2\pi^2G}\right)^l\\
 & \leq & \sum_{l=1}^{|J|}\left(\frac{d}{2\pi^2G}\right)^l\\
 & \leq & CdG^{-1}
\end{eqnarray*}
for some absolute constant $C>0$. Therefore we have
\begin{equation*}
 \sum_{||j||_{\infty}>G}\tilde{a}_j^2+\tilde{b}_j^2\leq CdG^{-1}
\end{equation*}
and the Lemma is proved.

\vspace{3ex}

With $L(N,G,\nu)$ as defined in (\ref{22}) we have

\begin{lemma}
\label{102}
 Let $f:\mathbb{R}^d\to\mathbb{R}$ be a periodic function satisfying $V_{HK}(f)\leq 1$ and let $(M_n)_{n\geq 0}$ be a lacunary sequence of matrices satisfying the Hadamard gap condition (\ref{17}). Then we have
 \begin{equation}
  \label{15}
  \mathlarger{\mathlarger{\int}}_{[0,1)^d}\left(\sum_{n=1}^Nf(M_nx)\right)^2\,dx\leq C(\log(d)||f||_2^2+||f||_2)N
 \end{equation}
 where $C>0$ is an absolute constant depending only on $q$. If the sequence $(M_n)_{n\geq 1}$ furthermore satisfies $L(N,G,0)=o(N)$
 for any fixed $G\geq 1$ then we have
 \begin{equation}
 \label{96}
 \lim_{N\to\infty}\frac{1}{N}\mathlarger{\mathlarger{\int}}_{[0,1)^d}\left(\sum_{n=1}^Nf(M_nx)\right)^2\,dx=||f||_2^2.
 \end{equation}
\end{lemma}

\vspace{3ex}

Note that hereafter we write $\log(x)$ for $\max(1,\log(x))$.

\vspace{3ex}

\textit{Proof.} For $||j||_{\infty},||j'||_{\infty}\leq G$ and $k>\log_q(G)$ we have
\begin{equation}
  \label{93}
 ||M^T_nj'||_{\infty}\leq ||M^T_n||_{\infty}||j'||_{\infty}\leq q^{-k}||M^T_{n+k}j||_{\infty}||j'||_{\infty}< ||M^T_{n+k}j||_{\infty}.
\end{equation}
Therefore we obtain $M^T_nj'\neq M^T_{n+k}j$. Now let $p_G$ be the $G$th Fej\'er mean of $f$. Then for some $k>\log_q(G)$ we have
\begin{equation*}
 p_G(M_nx)p_G(M_{n+k}x)=\sum_u\alpha_u\cos(2\pi\langle u,x\rangle)+\beta_u\sin(2\pi\langle u,x\rangle)
\end{equation*}
where any $u$ is of the form $M^T_nj\pm M^T_{n+k}j'$ for some $1\leq ||j||_{\infty},||j'||_{\infty}\leq G$. Therefore by (\ref{93}) we get
\begin{equation}
  \label{95}
 \int_{[0,1)^d}p_G(M_nx)p_G(M_{n+k}x)\,dx=0
\end{equation}
for $k>\log_q(G)$. By Lemma \ref{94} for any $k\geq 1$ there is a trigonometric polynomial $g_k$ with $dq^{2k}-1<\deg(g_k)\leq dq^{2k}$ such that
\begin{equation*}
 ||f-g_k||_2\leq Cq^{-k}.
\end{equation*}
Therefore for $k'>\log_q(dq^{2k})$ by (\ref{95}) and Cauchy-Schwarz inequality we have
\begin{equation}
\label{18}
\begin{aligned}
 \left|\int_{[0,1)^d} f(M_{n}x)f(M_{n+k'}x)\,dx\right| \:\:\: \leq & \:\:\: \left|\int_{[0,1)^d} (f-g_k)(M_{n}x)f(M_{n+k'}x)\,dx\right|\\
 &\:\:\: +\left|\int_{[0,1)^d} g_k(M_{n}x)g_k(M_{n+k'}x)\,dx\right|\\
 &\:\:\: +\left|\int_{[0,1)^d} g_k(M_{n}x)(f-g_k)(M_{n+k'}x)\,dx\right|\\
 \leq &\:\:\: 2C||f||_2q^{-k}
 \end{aligned}
\end{equation}
since $||g_k||_2\leq ||f||_2$. We obtain
\begin{multline*}
\mathlarger{\mathlarger{\int}}_{[0,1)^d}\left(\sum_{n=1}^Nf(M_nx)\right)^2\,dx\\
\begin{aligned}
\leq & N||f||_2^2+2\sum_{n=1}^N\sum_{k'=1}^{N-n}\left|\int_{[0,1)^d}f(M_nx)f(M_{n+k'}x)\,dx\right|\\
\leq & N||f||_2^2+2N(\log_q(d)+2)||f||_2^2+2\sum_{n=1}^N\sum_{k'=\log_q(d)+3}^{N-n}\left|\int_{[0,1)^d}f(M_nx)f(M_{n+k'}x)\,dx\right|\\
\leq & N||f||_2^2+CN\log_q(d)||f||_2^2+\sum_{n=1}^N\sum_{k=1}^{\infty}C||f||_2q^{-k}\\
\leq & C(\log(d)||f||_2^2+||f||_2)N.
\end{aligned}
\end{multline*}
Thus (\ref{15}) is shown.\\
The proof of (\ref{96}) is similar. For $k'>\log_q(q^{2k})$ instead of (\ref{18}) we get
\begin{equation*}
 \left|\int_{[0,1)^d} f(M_{n}x)f(M_{n+k}x)\,dx\right| \leq C||f||_2d^{1/2}q^{-k}
\end{equation*}
and similarly
\begin{equation*}
 \left|\int_{[0,1)^d} s_1(M_{n}x)s_2(M_{n+k}x)\,dx\right| \leq C||f||_2d^{1/2}q^{-k}
\end{equation*}
where for $i\in\{1,2\}$ the function $s_i$ is of the form $p_{G_i}$ or $r_{G_i}$ for some suitable number $G_i>0$.\\
For any $G\geq 1$ we obtain
\begin{multline*}
 \left|\frac{1}{N}\mathlarger{\mathlarger{\int}}_{[0,1)^d}\left(\sum_{n=1}^Nf(M_nx)\right)^2\,dx-||f||_2^2\right|\\
 \begin{aligned}
 \leq & \frac{2}{N}\sum_{n=1}^N\sum_{k=1}^{N-n}\left|\int_{[0,1)^d}f(M_nx)f(M_{n+k}x)\,dx\right|\\
 \leq & \frac{C}{N}\sum_{n=1}^N\sum_{k=1}^{N-n}||f||_2d^{1/2}\min(G^{-1/2},q^{-k})\\
 & +\frac{2}{N}\sum_{n=1}^N\sum_{k=1}^{N-n}\left|\int_{[0,1)^d}p_G(M_nx)p_G(M_{n+k}x)\,dx\right|\\
 \leq & C||f||_2d^{1/2}G^{-1/2}\\
 & + (2G+1)^{2d}\frac{h(N)}{N}
 \end{aligned}
\end{multline*}
for some function $h(N)$ with $h(N)/N\to 0$. Observe that such a function exists by assumption on $L(N,G,0)$. Since the constant $G\geq 1$ can be chosen arbitrary, (\ref{96}) is shown.

\begin{lemma}
 \label{135}
 Let $(M_n)_{n\geq 1}$ be some lacunary sequence satisfying the Hadamard gap condition (\ref{17}). For any $n,n'\in\mathbb{N}$ and $j\in\mathbb{Z}^d$ with $||j||_{\infty}\leq G$ for some $G\geq 2$ there exists at most one $j'\in\mathbb{Z}^d$ with $||j'||_{\infty}\leq G$ such that
 \begin{equation}
  \label{136}
  ||M^T_nj\pm M^T_{n'}j'||_{\infty}<||M^T_{n''}||_{\infty}
 \end{equation}
 where $n''\leq \min(n,n')-\log_q(G)$.
\end{lemma}

Here (\ref{136}) has to be understood in the following sense: For $n,n',j$ there exists at most one $j'$ with $||M^T_nj+M^T_{n'}j'||_{\infty}<||M^T_{n''}||_{\infty}$ and at most one possibly different
$j'$ with $||M^T_nj-M^T_{n'}j'||_{\infty}<||M^T_{n''}||_{\infty}$.

\vspace{3ex}

\textit{Proof.} Suppose that (\ref{136}) is satisfied for some $n,n',j,j'$. Now take some $j''\neq j'$. We observe
\begin{eqnarray*}
 ||M^T_nj\pm M^T_{n'}j''||_{\infty} & \geq & ||M^T_{n'}(j''-j')||_{\infty}-||M^T_nj\pm M^T_{n'}j'||_{\infty}\\
 & > & G||M^T_{n''}||_{\infty}-||M^T_{n''}||_{\infty}\\
 & \geq & ||M^T_{n''}||_{\infty}.
\end{eqnarray*}
Therefore the Lemma is proved.
 
\begin{lemma}
\label{123}
 Let $p(x)=\sum_{j\in\mathbb{Z}^d,0<||j||_{\infty}\leq G}a_j\cos(2\pi\langle j,x\rangle)+b_j\sin(2\pi\langle j,x\rangle)$ with $G\geq 2$ be some trigonometric polynomial satisfying $||p||_{\infty}\leq 1$ and $V_{HK}(p)\leq 1$. Let $(M_n)_{n\geq 1}$ be a lacunary sequence satisfying the Hadamard gap condition (\ref{17}). 
 Then we have
 \begin{equation}
  \label{108}
  \mathlarger{\mathlarger{\int}}_{[0,1)^d}\left(\sum_{n=1}^Np(M_nx)\right)^4\,dx\leq CN^2  
 \end{equation}
 for some constant $C>0$ which depends on $q$ and $d$. If furthermore $N\geq \max(Gd^{-1},d)$
then we have
 \begin{equation}
  \label{109}
  \mathlarger{\mathlarger{\int}}_{[0,1)^d}\left(\sum_{n=1}^Np(M_nx)\right)^4\,dx\leq Cd^{1/3}N^2
 \end{equation}
for some absolute constant $C>0$ depending only on $q$.    
\end{lemma}

\textit{Proof.} The Proof is this Lemma is based on \cite{BP79} where a similar result was proved for $d=1$.
We only prove (\ref{109}). The proof of (\ref{108}) is essentially the same.

First we are going to show
\begin{equation}
 \label{98}
 \mathbb{P}\left(\left|\sum_{n=1}^Np(M_nx)\right|>t\sqrt{N}\right)\leq C\exp\left(-\frac{t^{3/2}}{8C\log(d)}\right)
\end{equation}
for some constant $C>0$ depending only on $q$. For some $0<\beta<1$ set $P=\lfloor N^{\beta}\rfloor$ and $l=\lceil N/2P\rceil$.  Without loss of generality we may assume $3N^4\leq q^{N^{\beta}-1}$. There exists some $N_0\in\mathbb{N}$ depending only on $q$ with $3N^4\leq q^{N^{\beta}-1}$ for $N\geq N_0$. Therefore there is some constant $C_q>0$ which depends only on $q$ and $\beta$ such that (\ref{109}) is satisfied with $C_q$ and any $N<N_0$. Since $N\geq G^2$ we have $\log_q(G)+\log_q(3G)\leq \lfloor N^{\beta}\rfloor$ for $N\geq N_0$.
Now define
\begin{equation*}
 U_{m}(x)=\sum_{n=Pm+1}^{\min(P(m+1),N)}p(M_nx).
\end{equation*}
By Markov's inequality and Cauchy-Schwarz inequality we obtain
\begin{equation}
 \label{104}
 \begin{aligned}
 \mathbb{P}\left(\left|\sum_{n=1}^Np(M_nx)\right|>t\sqrt{N}\right) \:\:\: & \leq \:\:\: 2\exp(-\kappa_Nt\sqrt{N})\mathlarger{\mathlarger{\int}}_{[0,1)^d}\exp\left(\kappa_N\sum_{n=1}^Np(M_nx)\right)\,dx\\
 & \leq \:\:\: 2\exp(-\kappa_Nt\sqrt{N})I(\kappa_N,l)^{1/2}I'(\kappa_N,l)^{1/2} 
 \end{aligned}
\end{equation}
where $\kappa_N>0$ and
\begin{eqnarray*}
 I(\kappa_N,l) & = & \mathlarger{\mathlarger{\int}}_{[0,1)^d}\prod_{m=0}^{l-1}\exp\left(2\kappa_NU_{2m}(x)\right)\,dx,\\
 I'(\kappa_N,l) & = & \mathlarger{\mathlarger{\int}}_{[0,1)^d}\prod_{m=1}^{l}\exp\left(2\kappa_NU_{2m-1}(x)\right)\,dx.
\end{eqnarray*}
Since $I(\kappa_N,l)$ and $I'(\kappa_N,l)$ can be estimated similarly, we only estimate the first one.

Using $e^{|z|}\leq (1+z+z^2)e^{|z|^3}$ we observe
\begin{equation}
 \label{101}
 I(\kappa_N,l)\leq \exp\left(\sum_{m=0}^{l-1}|\kappa_NU_{2m}(x)|^3\right)\mathlarger{\mathlarger{\int}}_{[0,1)^d}\prod_{m=0}^{l-1}(1+\kappa_NU_{2m}(x)+(\kappa_NU_{2m}(x))^2)\,dx
\end{equation}
For some $m\geq 0$ we have
\begin{equation*}
 U_{2m}^2(x)=\sum_{j\in\mathbb{Z}^d}\alpha_j\cos(2\pi\langle j,x\rangle)+\beta_j\sin(2\pi\langle j,x\rangle)
\end{equation*}
for suitable numbers $\alpha_j$, $\beta_j$ for all $j\in\mathbb{Z}^d$. Now set
\begin{equation*}
 V_{2m}(x)=\sum_{j\in\mathbb{Z}^d,||j||_{\infty}< ||M_{m'}^T||_{\infty}}\alpha_j\cos(2\pi\langle j,x\rangle)+\beta_j\sin(2\pi\langle j,x\rangle)
\end{equation*}
for $m'=\lfloor 2Pm-\log_q(G)\rfloor$ and furthermore $W_{2m}(x)=U_{2m}^2(x)-V_{2m}(x)$. To estimate $|V_{2m}(x)|$ observe that by Lemma \ref{135} for any $1\leq n,n'\leq N$ and $j\in\mathbb{Z}^d$ with $||j||_{\infty}\leq G$ there is at most one $j'\in\mathbb{Z}^d$ with $||j'||_{\infty}\leq G$ such that $||M_n^Tj\pm M_{n'}^Tj'||_{\infty}<||M_{n''}^T||_{\infty}$ where $n''\leq \min(n,n')-\log_q(G)$.
For $||j||_{\infty},||j'||_{\infty}\leq 1/2\cdot q^{|n-n'|}$ we have
\begin{equation*}
 ||M^T_nj\pm M^T_{n'}j'||_{\infty}\geq\frac{1}{2}q^{\min(n,n')-n''}||M^T_{n''}||_{\infty} \geq ||M^T_{m'}||_{\infty}.
\end{equation*}
Then by Lemma \ref{40} and Cauchy-Schwarz inequality we have
\begin{equation}
 \label{142}
 |V_{2m}(x)|\leq\sum_{n,n'=1}^N\sum_{\substack {1\leq ||j||_{\infty},||j'||_{\infty}\leq G,\\ ||M^T_nj\pm M^T_{n'}j'||_{\infty}<||M^T_{m'}||_{\infty}}}(|a_j|+|b_j|)(|a_{j'}|+|b_{j'}|)\leq C\sqrt{d}||p||_2N
\end{equation}
where the constant $C>0$ depends only on $q$.
Since
\begin{equation*}
 0\leq 1+\kappa_nU_{2m}(x)+(\kappa_nU_{2m}(x))^2
\end{equation*}
we obtain
\begin{multline}
 \label{103}
 I(\kappa_N,l)\leq \exp(C\kappa_N^3N^{1+2\beta})\\
 \cdot\mathlarger{\mathlarger{\int}}_{[0,1)^d}\prod_{m=0}^{l-1}(1+C\sqrt{d}||p||_2)P+\kappa_NU_{2m}(x)+\kappa_N^2W_{2m}(x))\,dx.
\end{multline}
Furthermore we get
\begin{equation*}
 \kappa_NU_{2m}(x)+\kappa_N^2W_{2m}(x)=\sum_{j}\alpha_j\cos(2\pi\langle j,x\rangle)+\beta_j\sin(2\pi\langle j,x\rangle)
\end{equation*}
for suitable $\alpha_j,\beta_j$ and $j\in\mathbb{Z}^d$ satisfying $||M_{m'}^T||_{\infty}\leq ||j||_{\infty}\leq 2G||M_{P(2m+1)}^T||_{\infty}$.

Now for $0\leq k\leq m$ let $j_{2k}$ be any frequency vector of the trigonometric polynomial $1+C(\log(d)||p||_2^2+||p||_2)N+\kappa_NU_{2m}(x)+\kappa_N^2W_{2m}(x)$. If $j_{2m}\neq 0$ then we get
\begin{multline}
 \label{131}
  ||j_{2m}||_{\infty}-\sum_{k=0}^{m-1}||j_{2k}||_{\infty}\\
 \begin{aligned}
 & \geq ||M_{m'}^T||_{\infty}-2G\sum_{k=0}^{m-1}||M^T_{P(2k+1)||_{\infty}}\\
 & \geq 3G||M^T_{P(2m-1)}||_{\infty}-2G||M^T_{P(2m-1)}||_{\infty}\sum_{k=0}^{m-1}q^{P(2k+1)-P(2m-1)}\\
 & \geq \left(3-2\frac{1}{1-q^{-2P}}\right)G||M_{P(2m-1)}^T||_{\infty}\\
 & > 0
 \end{aligned}
\end{multline}
where the second inequality follows by assumption on $P>\log_q(G)+\log_q(3G)$ and the fourth inequality follows by $P>\log_q(\sqrt{3})$ which without loss of generality we may assume. Therefore by (\ref{103}) we have
\begin{equation*}
 I(\kappa_N,l)\leq \exp(C\kappa_N^3N^{1+2\beta})\mathlarger{\mathlarger{\int}}_{[0,1)^d}\prod_{m=0}^{l-1}(1+C\sqrt{d}||p||_2P)\,dx.
\end{equation*}
Plugging this into (\ref{104}) yields
\begin{equation}
 \label{106}
 \mathbb{P}\left(\left|\sum_{n=1}^Np(M_nx)\right|>t\sqrt{N}\right)\leq 2\exp\left(-\kappa_Nt\sqrt{N}+C_0\kappa_N^3N^{1+2\beta}+C_0\sqrt{d}\kappa_N^2N\right)
\end{equation}
for some absolute constant $C_0$. Choose $\beta=1/4$. Then for 
\begin{equation*}
 \kappa_N=\frac{t}{2C_0\sqrt{d}\sqrt{N}}
\end{equation*}
and $0\leq t\leq C_0d$ we observe
\begin{multline}
 \label{105}
 \exp\left(-\kappa_Nt\sqrt{N}+C_0\kappa_N^3N^{1+2\beta}+C_0\sqrt{d}\kappa_N^2N\right)\\
 \begin{aligned}
  & = \exp\left(C_0\sqrt{d}\frac{t^2}{4C_0^2d}+C_0\frac{t^3}{8C_0^3d^{3/2}}-\frac{t^2}{2C_0\sqrt{d}}\right)\\
  & = \exp\left(-\frac{t^2}{4C_0\sqrt{d}}\left(1-\frac{t}{2C_0d}\right)\right)\\
  & \leq \exp\left(-\frac{t^2}{8C_0\sqrt{d}}\right).
 \end{aligned}
\end{multline}
A similar calculation for 
\begin{equation*}
 \kappa_N=\frac{\sqrt{t}}{2\sqrt{C_0N}}
\end{equation*} 
and $t>C_0d$ yields
\begin{multline}
 \label{107}
 \exp\left(-\kappa_Nt\sqrt{N}+C_0\kappa_N^3N^{1+2\beta}+C_0\sqrt{d}\kappa_N^2N\right)\\
 \begin{aligned}
  & = \exp\left(C_0\sqrt{d}\frac{t}{4C_0}+C_0\frac{t^{3/2}}{8C_0^{3/2}}-\frac{t^{3/2}}{2C_0^{1/2}}\right)\\
  & = \exp\left(-\frac{t^{3/2}}{8C_0^{1/2}}\left(3-\frac{2C_0^{1/2}}{t^{1/2}}\sqrt{d}\right)\right)\\
  & \leq \exp\left(-\frac{t^{3/2}}{8C_0^{1/2}}\right).
 \end{aligned}
\end{multline}
Combining (\ref{106}), (\ref{105}) and (\ref{107}) we get (\ref{98}). Therefore we have
\begin{eqnarray*}
 \mathlarger{\mathlarger{\int}}_{[0,1)^d}\left(\sum_{n=1}^Np(M_nx)\right)^4\,dx & \leq & \mathlarger{\mathlarger{\int}}_{0}^{\infty}tN^2d\mathbb{P}\left(\left|\sum_{n=1}^Np(M_nx)\right|>t^{1/4}N^{1/2}\right)\\
  & \leq & N^2\mathlarger{\mathlarger{\int}}_{0}^{\infty}Ct\frac{t^{1/2}}{\sqrt{d}}\exp\left(-\frac{t^{3/2}}{8C\sqrt{d}}\right)\,dt\\
  & \leq &  N^2\mathlarger{\mathlarger{\int}}_{0}^{\infty}C(t')^{3/2}\exp\left(-\frac{(t')^{3/2}}{8C}\right)d^{1/3}\,dt'
\end{eqnarray*}
where in the last line we substituted $t'=t/(d)^{1/3}$. Thus we observe
\begin{equation*}
 \mathlarger{\mathlarger{\int}}_{[0,1)^d}\left(\sum_{n=1}^Np(M_nx)\right)^4\,dx \leq Cd^{1/3}N^2
\end{equation*}
for some constant $C>0$ which only depends on $q$. Therefore the Lemma is proved.


\vspace{5ex}




\vspace{5ex}

\section{Proof of Theorem (\ref{3})}
\label{S31}

The proof of Theorem \ref{3} is essentially based on the following Theorem due to Heyde and Brown \cite{HB70} which is a consequence of Strassen's almost sure invariance principle for martingale differences sequences. We will use a generalized version stated in \cite{AB10}.
\begin{theorem}[\extcite{AB10}{Theorem B},\cite{HB70}]
\label{5}
 Let $(X_k,\mathcal{F}_k,k\geq 1)$ be a martingale differences sequence with finite fourth moments. Set
$V_K=\sum_{k=1}^K\mathbb{E}[X_k^2|\mathcal{F}_{k-1}]$ and let $(b_K)_{K\geq 1}$ be a sequence of positive numbers. Then we have
\begin{equation}
\label{14}
 \sup_t\left|\mathbb{P}\left(\frac{1}{\sqrt{b_K}}\sum_{k=1}^KX_k<t\right)-\Phi(t)\right|\leq A\left(\frac{\sum_{k=1}^K\mathbb{E}[X_k^4]+\mathbb{E}\left[(V_K-b_K)^2\right]}{b_K^2}\right)^{1/5},
\end{equation}
where $A$ is an absolute constant.
\end{theorem}

\textit{Proof of Theorem \ref{3}.} First we are going to show (\ref{139}). Therefore for some fixed but large enough integer $N\geq 1$ set $G=\lfloor d\max(2,N^{\alpha})\rfloor$ for some $0<\alpha\leq 1$ and define $p=p_G$ and $r=r_G$ as in (\ref{92}). Without loss of generality we may assume $||f||_{\infty}\leq 1$ and $V_{HK}(f)\leq 1$. Therefore it is easy to see that $||p||_2\leq ||f||_2\leq 1$ and $||p||_{\infty}\leq ||f||_{\infty}\leq 1$ as well as $||r||_{\infty}\leq ||f||_{\infty}+||p||_{\infty}\leq 2$.
We now decompose the set $\{1,\ldots,N\}$ into consecutive blocks $\Delta_1,\Delta'_2,\Delta_2,\ldots,\Delta'_k,$ $\Delta_k,\ldots$ such that the blocks have length
\begin{equation*}
 |\Delta'_k|=\lceil 2(1+2\eta)\log_q(k)+C_{d,G}\rceil, \quad |\Delta_k|=\lfloor12\eta^{-1}C_{d,G}^{1+\eta}k^{\eta}\rfloor.
\end{equation*}
for some $0<\eta<1$. The constant $C_{d,G}$ is defined by
\begin{equation}
\label{110}
C_{d,G}=C(\log(G)+\log(d)+d\log(5\log(G)))
\end{equation}
for some large enough $C$ depending only on $q$. It can easily be shown that $|\Delta_{k+1}'|\leq |\Delta_k|$ for all $k\geq 1$.
In order to define a suitable martingale differences sequence we replace $f$ by its low-frequency part $p$ which is a finite trigonometric polynomial. Furthermore we will neglect the indices in $\Delta'_k$. The purpose of this is having a fast enough decreasing ratio
\begin{equation*}
 \frac{||M^T_{(k-1)^+}||_{\infty}}{||M^T_{k^-}||_{\infty}}\leq k^{-2(1+2\eta)}q^{-C_{d,G}},
\end{equation*}
where $k^+,k^-$ is the largest resp. the smallest integer in $\Delta_k$. Later on it will be shown that the asymptotic size of $\sum_{n=1}^Nr(M_nx)$ and $\sum_k\sum_{n\in\Delta'_k}p(M_nx)$ can be neglected. We now approximate $p(M_nx)$ by a piecewise constant function $\varphi_n(x)$ which is necessary to define the martingale differences sequence.
Let 
\begin{equation*}
 m(n)=\lceil \log_2(||M^T_{k^+}||_{\infty})+(1+2\eta)\log_2(k)+C'_{d,G}\rceil
\end{equation*}
where $k=k_n$ is defined by $n\in\Delta_{k_n}$ and $C'_{d,G}$ is a constant depending on $q$, $d$ and $G$ such that 
\begin{multline}
\label{111}
 2(1+\eta)\log_2(C_{d,G})+\log_2(G)+\log_2(d)+(2d+1)\log_2(5\log(G))\\
 \begin{aligned}
 & \leq C_{d,G}'\\
 & \leq \log_2(q)C_{d,G}-2(1+\eta)\log_2(C_{d,G})-d\log_2(5\log(G)).
 \end{aligned}
\end{multline}
Observe that such a constant $C_{d,G}'$ exists if $C_{d,G}$ in (\ref{110}) is chosen large enough.
Let $\mathcal{F}_k$ be the $\sigma$-field generated all sets of the form
\begin{equation*}
 \left[\frac{v_1}{2^{m(k^+)}},\frac{v_1+1}{2^{m(k^+)}}\right)\times\cdots\times\left[\frac{v_d}{2^{m(k^+)}},\frac{v_d+1}{2^{m(k^+)}}\right)
\end{equation*}
with $v_i\in\{0,1,\ldots,2^{m(k^+)}-1\}$ for $i\in\{1,\ldots,d\}$.
With $x,x'\in\mathcal{A}_k$ for any atom $\mathcal{A}_k=\prod_{i=1}^d\mathcal{A}_{k,i}\in\mathcal{F}_k$ we have
\begin{eqnarray*}
 |p(M_nx)-p(M_nx')| & \leq & \sum_{i=1}^d \sup_{y\in\mathcal{A}_k}\left|\frac{\partial}{\partial y_i}p(M_ny)\right|2^{-m(n)}\\
  & \leq & Cd||M_n^T||_{\infty}G(5\log(G))^d2^{-m(n)}\\
  & \leq & C\cdot C_{d,G}^{-2(1+\eta)}k^{-(1+2\eta)}.
\end{eqnarray*}
Here the second inequality follows by Lemma \ref{40}. The third inequality follows by (\ref{111}).
Therefore on any atom of $\mathcal{F}_k$ we easily can find some constant function, say $\hat{\varphi_n}$, such that
\begin{equation}
 \label{113}
 |p(M_nx)-\hat{\varphi}_n(x)|\leq C\cdot C_{d,G}^{-2(1+\eta)}k^{-(1+2\eta)}.
\end{equation}
We have
\begin{equation*}
 p(M_nx)=\sum_{1\leq ||j||_{\infty}\leq G}a_j'\cos(2\pi\langle M_n^Tj,x\rangle)+b_j'\sin(2\pi\langle M_n^Tj,x\rangle).
\end{equation*}
Thus for any atom $\mathcal{A}_{k-1}\subset \mathcal{F}_{k-1}$ we obtain
\begin{multline*}
 \frac{1}{\lambda(\mathcal{A}_{k-1})}\left|\int_{\mathcal{A}_{k-1}}p(M_nx)\,dx\right|\\
 \leq\frac{1}{\lambda(\mathcal{A}_{k-1})}\sum_{1\leq ||j||_{\infty}\leq G}|a_j'|\left|\int_{\mathcal{A}_{k-1}}\cos(2\pi\langle M_n^Tj,x\rangle)\,dx\right|+|b_j'|\left|\int_{\mathcal{A}_{k-1}}\sin(2\pi\langle M_n^Tj,x\rangle)\,dx\right|
\end{multline*}
where $\lambda$ denotes the Lebesgue-measure on $[0,1)^d$. For fixed $n$ and $j$ choose $i\in\{1,\ldots,d\}$ such that $|(M_n^Tj)_i|=||M^T_nj||_{\infty}$. Then we observe by angle sum and difference identities of trigonometric functions that
\begin{equation*}
\int_{B_k}\cos(2\pi\langle M_n^Tj,x\rangle)\,dx=\int_{B_k}\sin(2\pi\langle M_n^Tj,x\rangle)\,dx=0
\end{equation*}
for any box $B_k$ such that $B_{k,i'}=\mathcal{A}_{k,i'}$ for $i\neq i'$ and $\lambda(B_{k,i})=v/||M_n^Tj||_{\infty}$ for some integer $v$. Thus we have
\begin{multline}
 \label{146}
 \frac{1}{\lambda(\mathcal{A}_{k-1})}\left|\int_{\mathcal{A}_{k-1}}p(M_nx)\,dx\right|\\
 \begin{aligned}
 & \leq \sum_{1\leq ||j||_{\infty}\leq G}\frac{|a_j'|+|b_j'|}{||M_n^Tj||_{\infty}}2^{m((k-1)^+)}\\
 & \leq \sum_{1\leq ||j||_{\infty}\leq G}(|a_j'|+|b_j'|)\cdot \frac{||M_{(k-1)^+}^T||_{\infty}}{q^{C_{d,G}}||M_{k^--C_{d,G}}^T||_{\infty}}k^{1+2\eta}2^{C'_{d,G}}\\
 & \leq C\cdot (5\log(G))^dq^{-C_{d,G}}k^{-(1+2\eta)}2^{C'_{d,G}}\\
 & \leq C\cdot C_{d,G}^{-2(1+\eta)}k^{-(1+2\eta)}.
 \end{aligned}
\end{multline}
Now set
\begin{equation*}
 \varphi_n(x)=\hat{\varphi}_n(x)-\mathbb{E}[\hat{\varphi}_n|\mathcal{F}_{k-1}].
\end{equation*}
Then by (\ref{113}) we have
\begin{equation}
\label{115}
 \begin{aligned}
 |p(M_nx)-\varphi_n(x)| \:\:\: & \leq \:\:\: |p(M_nx)-\hat{\varphi}_n(x)|+|\hat{\varphi}_n(x)-\varphi_n(x)|\\
 & \leq \:\:\: C\cdot C_{d,G}^{-2(1+\eta)}k^{-(1+2\eta)}
 \end{aligned}
\end{equation}
for some absolute constant $C>0$ depending only on $q$.
Thus $p(M_nx)$ can be approximated by a function $\varphi_n(x)$, which is constant on any atom of $\mathcal{F}_k$ and satisfies:
\begin{enumerate}
 \item $|p(M_nx)-\varphi_n(x)|\leq C\cdot C_{d,G}^{-2(1+\eta)}k^{-(1+2\eta)} \quad \textnormal{for all } x\in[0,1)^d,$
 \item $\mathbb{E}[\varphi_n|\mathcal{F}_{k-1}]=0 \quad \textnormal{for all } n\in\Delta_k.$
\end{enumerate}
We now define
\begin{equation*}
 X_k=\sum_{n\in\Delta_k}\varphi_n(x), \quad V_{K+1}=\sum_{k=1}^{K+1}\mathbb{E}[X_k^2|\mathcal{F}_{k-1}],
\end{equation*}
where $K$ is given such that $N\in\Delta_{K+1}'\cup\Delta_{K+1}$. It is easy to see that $(X_k,\mathcal{F}_k)$ is a martingale differences sequence. We are going to approximate $\sum_{n=1}^Nf(M_nx)$ by $\sum_{k=1}^KX_k$ in order to apply Theorem \ref{5}.
Therefore we need some more notations. For $k\in\{1,\ldots,K+1\}$ define
\begin{equation*}
 Y_k=\sum_{n\in\Delta_k}p(M_nx), \quad Y'_k=\sum_{n\in\Delta'_k}p(M_nx).
\end{equation*}
We observe
\begin{equation*}
N\leq \sum_{k=1}^{K+1}|\Delta'_k|+|\Delta_k|\leq 2\sum_{k=1}^{K+1}|\Delta_k|\leq C\eta^{-1}\cdot C_{d,G}^{1+\eta}K^{1+\eta}
\end{equation*}
and
\begin{equation*}
N\geq \sum_{k=1}^{K}|\Delta'_k|+|\Delta_k|\geq \sum_{k=1}^{K}|\Delta_k|\geq C\eta^{-1}\cdot C_{d,G}^{1+\eta}K^{1+\eta}.
\end{equation*}
There are constants $C_1,C_2>0$ such that
\begin{equation}
 \label{114}
 C_1\eta^{-1/(1+\eta)}C_{d,G}^{-1}N^{1/(1+\eta)}\leq K \leq C_2\eta^{-1/(1+\eta)}C_{d,G}^{-1}N^{1/(1+\eta)}.
\end{equation}

By definition we have the following decomposition
\begin{equation*}
 \sum_{n=1}^Nf(M_nx)=\sum_{k=1}^{K+1} X_k+\sum_{k=1}^{K+1}(Y_k-X_k)+\sum_{k=1}^{K+1}Y_k'+\sum_{n=1}^Nr(M_nx).
\end{equation*}
Then standard estimates give
\begin{equation}
\label{50}
 \begin{aligned}
  \mathbb{P}\left(\sum_{n=1}^Nf(M_nx)\leq t\sigma_N\right) \:\:\: \leq & \:\:\: \mathbb{P}\left(\sum_{k=1}^{K+1} X_k\leq (t+\varepsilon)\sigma_N\right)+\mathbb{P}\left(\left|\sum_{k=1}^{K+1}(X_k-Y_k)\right|>\frac{\varepsilon\sigma_N}{3}\right)\\
  & \:\:\: +\mathbb{P}\left(\left|\sum_{k=1}^{K+1} Y'_k\right|>\frac{\varepsilon\sigma_N}{3}\right)+\mathbb{P}\left(\left|\sum_{n=1}^Nr(M_nx)\right|>\frac{\varepsilon\sigma_N}{3}\right),
 \end{aligned}
\end{equation}
for some $\varepsilon>0$ as well as
\begin{equation*}
 \begin{aligned}
  \mathbb{P}\left(\sum_{n=1}^Nf(M_nx)\geq t\sigma_N\right) \:\:\: \geq & \:\:\: \mathbb{P}\left(\sum_{k=1}^{K+1} X_k\geq (t-\varepsilon)\sigma_N\right)-\mathbb{P}\left(\left|\sum_{k=1}^{K+1}(X_k-Y_k)\right|>\frac{\varepsilon\sigma_N}{3}\right)\\
  & \:\:\: -\mathbb{P}\left(\left|\sum_{k=1}^{K+1} Y'_k\right|>\frac{\varepsilon\sigma_N}{3}\right)-\mathbb{P}\left(\left|\sum_{n=1}^Nr(M_nx)\right|>\frac{\varepsilon\sigma_N}{3}\right).
 \end{aligned}
\end{equation*}
Since both inequalities can be estimated analogously we will only focus on the first one. First we estimate the three latter terms before we estimate the first one by applying Theorem \ref{5}.
In order to estimate the second term by (\ref{115}) we have
\begin{eqnarray*}
 \left|\sum_{k=1}^{K+1}(X_k-Y_k)\right| & = & \left|\sum_{k=1}^{K+1}\left(\sum_{n\in\Delta_k}\varphi_n(x)-\sum_{n\in\Delta_k}p(M_nx)\right)\right|\\
 & \leq & \sum_{k=1}^{K+1}C|\Delta_k|\cdot C_{d,G}^{-2(1+\eta)}k^{-(1+2\eta)}\\
 & \leq & C\eta^{-1}C_{d,G}^{-(1+\eta)}
\end{eqnarray*}
for some constant $C>0$ depending only on $q$.
Therefore by Chebyshev's inequality we see that
\begin{equation}
\label{46}
\mathbb{P}\left(\left|\sum_{k=1}^{K+1}(X_k-Y_k)\right|>\frac{\varepsilon\sigma_N}{3}\right)\leq C\eta^{-2}C_{d,G}^{-2(1+\eta)}\varepsilon^{-2}N^{-1}.
\end{equation}
To estimate the third term we use the definition of $Y_k'$, Lemma \ref{102} and (\ref{114}) to get
\begin{eqnarray*}
 \left|\left|\sum_{k=1}^{K+1} Y'_k\right|\right|_2^2 & \leq & C\sum_{k=1}^{K+1}|\Delta'_k|\log(d)\\
 & \leq  & C\log(d)N^{1/(1+\eta)}\log(N).
\end{eqnarray*}
Thus by (\ref{21}) and Chebyshev's inequality we obtain
\begin{equation}
\label{47}
 \mathbb{P}\left(\left|\sum_{k=1}^{K+1} Y'_k\right|>\frac{\varepsilon\sigma_N}{3}\right)\leq C\log(d)\varepsilon^{-2}N^{-\eta/(1+\eta)}\log(N).
\end{equation}
With another application of Chebyshev's inequality, Lemma \ref{94} and Lemma \ref{102} we get the following estimate for the fourth term
\begin{equation}
\label{48}
 \mathbb{P}\left(\left|\sum_{n=1}^Nr(M_nx)\right|>\frac{\varepsilon\sigma_N}{3}\right)\leq C\log(d)\varepsilon^{-2}d^{1/2}G^{-1/2}.
\end{equation}
Thus by (\ref{50}),(\ref{46}),(\ref{47}) and (\ref{48}) we have
\begin{multline}
 \label{116}
 \mathbb{P}\left(\sum_{n=1}^Nf(M_nx)\leq t\sigma_N\right)\\
 \leq\mathbb{P}\left(\sum_{k=1}^{K+1} X_k\leq (t+\varepsilon)\sigma_N\right)+C\eta^{-2}\log(d)\varepsilon^{-2}(N^{-\eta/(1+\eta)}\log(N)+d^{1/2}G^{-1/2}).
\end{multline}
Now for any integer $k\in\{1,\ldots,K+1\}$ define
\begin{equation*}
 s_{k}^2=\sum_{l=1}^{k}\mathlarger{\mathlarger{\int}}_{[0,1]^d}\left(\sum_{n\in\Delta_l}p(M_nx)\right)^2\, dx.
\end{equation*}
Thus by  we obtain
\begin{multline}
\label{52}
 \left|\mathbb{P}\left(\sum_{n=1}^Nf(M_nx)\leq t\sigma_N\right)-\Phi(t)\right|\\
  \begin{aligned}
  \leq & \left|\mathbb{P}\left(\sum_{k=1}^{K+1}X_k\leq s_{K+1}(t+\varepsilon)\frac{\sigma_N}{s_{K+1}}\right)-\Phi\left((t+\varepsilon)\frac{\sigma_N}{s_{K+1}}\right)\right|\\
  & +\left|\Phi\left((t+\varepsilon)\frac{\sigma_N}{s_{K+1}}\right)-\Phi(t)\right|\\
  & +C\eta^{-2}\log(d)\varepsilon^{-2}(N^{-\eta/(1+\eta)}\log(N)+d^{1/2}G^{-1/2}).
 \end{aligned}
\end{multline}

To apply Theorem \ref{5} we need a sequence of positive numbers $(b_K)_{K\geq 1}$. In \cite{HB70} the sequence was given by $\sum_{k=1}^K\mathbb{E}\left[\mathbb{E}[X_k^2|\mathcal{F}_{k-1}]\right]$. Here we take $s_K^2$ instead since later we are going to estimate the conditional second moments of $X_k$ by those of $Y_k$. 
In order to estimate the second term on the right-hand side of (\ref{52}) we now show that for some $C>0$ we have
\begin{equation}
\label{53}
 |\sigma_N-s_{K+1}|\leq (C\eta^{-1/(1+\eta)}\log(d)(N^{1/(1+\eta)}\log(N)+Nd^{1/2}G^{-1/2}))^{1/2}.
\end{equation}
Therefore we use standard estimates and observe
\begin{multline}
 \label{55}
 \left|s_{K+1}^2-\mathlarger{\mathlarger{\int}}_{[0,1)^d}\left(\sum_{n\in\cup_{k=1}^{K+1}}p(M_nx)\right)^2\, dx\right|\\
 \begin{aligned}
  & \leq \left|\mathlarger{\mathlarger{\int}}_{[0,1)^d}\left(\sum_{k=1}^{K+1}\left(\sum_{n\in \Delta_k}p(M_nx)\right)^2-\left(\sum_{n\in\cup_{k}\Delta_k}p(M_nx)\right)^2\right)\,dx\right|\\
  & \leq \sum_{\substack{n,n'\in\cup_k\Delta_k,\\ (n,n')\notin\Delta_k\times\Delta_k}}\left|\int_{[0,1)^d}p(M_nx)p(M_{n'}x)\,dx\right|\\
  & = 0
 \end{aligned}
\end{multline}
since $||M^T_nj\pm M^T_{n'}j'||_{\infty}>0$ if $1\leq ||j||_{\infty},||j'||_{\infty}\leq G$ and $|n-n'|\geq \log_q(G)$.
We now decompose
\begin{equation*}
 \sum_{n=1}^N f(M_nx)=\sum_{n\in\cup_{k=1}^{K+1}\Delta'_k}p(M_nx)+\sum_{n\in\cup_{k=1}^{K+1}\Delta_k}p(M_nx)+\sum_{n=1}^{N}r(M_nx)
\end{equation*}
By (\ref{55}) we have
\begin{equation}
\label{117}
 \left|\left|\sum_{n\in\cup_{k=1}^{K+1}\Delta_k}p(M_nx)\right|\right|_2\leq s_{K+1}.
\end{equation}
Now with Lemma \ref{102} and (\ref{114}) we obtain
\begin{equation}
 \label{118}
 \left|\left|\sum_{n\in\cup_{k=1}^{K+1}\Delta'_k}p(M_nx)\right|\right|_2\leq\left(C\eta^{-1/(1+\eta)}\log(d)N^{1/(1+\eta)}\log(N)\right)^{1/2}.
\end{equation}
By Lemma \ref{94} and \ref{102} we observe
\begin{equation}
 \label{119}
 \left|\left|\sum_{n=1}^Nr(M_nx)\right|\right|_2\leq\left(C\log(d)d^{1/2}G^{-1/2}N\right)^{1/2}.
\end{equation}
Therefore we obtain (\ref{53}). If we choose $N\in\mathbb{N}$ large enough such that
\begin{equation}
\label{121}
\frac{N}{N^{1/(1+\eta)}\log(N)+d^{1/2}G^{-1/2}N}\geq \frac{9C''\eta^{-1/(1+\eta)}\log(d)}{C'}
\end{equation}
where the constant $C',C''>0$ are defined such that (\ref{21}) is satisfied with $C=C'$ resp. (\ref{53}) is satisfied $C=C''$ then we obtain
\begin{equation*}
|\sigma_N-s_{K+1}|\leq \frac{\sqrt{C'N}}{3}.
\end{equation*}
Immediately we observe
\begin{equation}
\label{38}
s_{K+1}\geq \sigma_N-|\sigma_N-s_{K+1}|\geq\frac{2\sqrt{C'N}}{3}.
\end{equation} 
We furthermore get
\begin{equation}
\label{125}
\left|\frac{\sigma_N}{s_{K+1}}-1\right|\leq (C\eta^{-1/(1+\eta)}\log(d)(N^{-\eta/(1+\eta)}\log(N)+d^{1/2}G^{-1/2}))^{1/2}.
\end{equation}

Thus it can easily be shown that for a suitable constant $C>0$ depending only on $q$ and some large enough $N$ fulfilling (\ref{121}) by Mean Value Theorem we have
\begin{multline}
 \label{54}
 \left|\Phi\left((t+\varepsilon)\frac{\sigma_N}{s_{K+1}}\right)-\Phi(t)\right|\\
 \begin{aligned}
  & \leq \left(\frac{\sigma_N}{s_{K+1}}\varepsilon+\left|\left(\frac{\sigma_N}{s_{K+1}}-1\right)t\right|\right)\sup\left\{\frac{1}{\sqrt{2\pi}}e^{-1/2\cdot u^2}:t\leq u\leq (t+\varepsilon)\frac{\sigma_N}{s_{K+1}}\right\}\\
  & \leq (C\eta^{-1/(1+\eta)}\log(d))^{1/2}(\varepsilon+(N^{-\eta/(1+\eta)}\log(N)+d^{1/2}G^{-1/2})^{1/2}).
 \end{aligned}
\end{multline}
We plug this into (\ref{52}) and get
\begin{multline}
 \label{122}
 \left|\mathbb{P}\left(\sum_{n=1}^Nf(M_nx)\leq t\sigma_N\right)-\Phi(t)\right|\\
  \begin{aligned}
  \leq & \left|\mathbb{P}\left(\sum_{k=1}^{K+1}X_k\leq s_{K+1}t\right)-\Phi\left(t\right)\right|\\
  & +(C\eta^{-1/(1+\eta)}\log(d))^{1/2}(\varepsilon+(N^{-\eta/(1+\eta)}\log(N)+d^{1/2}G^{-1/2})^{1/2})\\
  & +C\eta^{-2}\log(d)\varepsilon^{-2}(N^{-\eta/(1+\eta)}\log(N)+d^{1/2}G^{-1/2}).
 \end{aligned}
\end{multline}
Therefore it remains to estimate
\begin{equation*}
 \left|\mathbb{P}\left(\sum_{k=1}^{K+1}X_k\leq s_{K+1}t\right)-\Phi\left(t\right)\right|
\end{equation*}
for which we use Theorem \ref{5}.
By Lemma \ref{123} we easily see
\begin{equation}
\label{13}
 \sum_{k=1}^{K+1}\mathbb{E}[X_k^4]\leq\sum_{k=1}^{K+1}C\log(d)^{2/3}|\Delta_k|^2\leq 
 C\eta^{-2}\log(d)^{2/3}C_{d,G}^{2(1+\eta)}K^{1+2\eta}.
\end{equation}
We define
\begin{equation*}
 \varsigma_k=\mathlarger{\mathlarger{\int}}_{[0,1)^d}\left(\sum_{n\in\Delta_k}p(M_nx)\right)^2\,dx.
\end{equation*}
By (\ref{113}) we observe $|X_k^2-Y_k^2|\leq C\eta^{-2}k^{-1}$ and therefore we obtain
\begin{equation}
\label{126} ||V_{K+1}-s_{K+1}^2||_2\leq\left|\left|\sum_{k=1}^{K+1}\mathbb{E}\left[Y_k^2-\varsigma_k|\mathcal{F}_{k-1}\right]\right|\right|_2+C\eta^{-2}\log(K).
\end{equation}
We now are going to decompose the terms $Y_k^2-\varsigma_k$. Therefore we set
\begin{eqnarray*}
\mathcal{R}_k(x) & = & \sum_{n\in\Delta_k}\sum_{\substack{n'\in\Delta_k,\\ |n-n'|\leq 1+\log_q(G+1)}}p(M_nx)p(M_{n'}x),\\
Q_k(x) & = & Y_k^2-\mathcal{R}_k(x).
\end{eqnarray*}
Furthermore we define
\begin{equation}
 \label{130}
 \begin{aligned}
  T^+(j,j',n,n',x) = & \frac{a'_{j}a'_{j'}}{2}\cos(2\pi\langle M_n^Tj+M_{n'}^Tj',x\rangle)+\frac{a'_{j}b'_{j'}}{2}\sin(2\pi\langle M_n^Tj+M_{n'}^Tj',x\rangle)\\
  & +\frac{b'_{j}a'_{j'}}{2}\sin(2\pi\langle M_n^Tj+M_{n'}^Tj',x\rangle)-\frac{b'_{j}b'_{j'}}{2}\cos(2\pi\langle M_n^Tj+M_{n'}^Tj',x\rangle),\\
  T^-(j,j',n,n',x) = & \frac{a'_{j}a'_{j'}}{2}\cos(2\pi\langle M_n^Tj-M_{n'}^Tj',x\rangle)+\frac{a'_{j}b'_{j'}}{2}\sin(2\pi\langle M_n^Tj-M_{n'}^Tj',x\rangle)\\
  & +\frac{b'_{j}a'_{j'}}{2}\sin(2\pi\langle M_n^Tj-M_{n'}^Tj',x\rangle)-\frac{b'_{j}b'_{j'}}{2}\cos(2\pi\langle M_n^Tj-M_{n'}^Tj',x\rangle),
\end{aligned}
\end{equation}
for $1\leq||j||_{\infty},||j'||_{\infty}\leq G$, $n,n'\in\mathbb{N}$ and $x\in [0,1)^d$.
We set
\begin{equation}
 \begin{aligned}
  \label{143}
  R_k(x) \:\:\: = & \:\:\: \sum_{\substack{n,n'\in\Delta_k,\\ |n-n'|\leq 1+\log_q(1+G)}}\sum_{\substack{1\leq||j||_{\infty},||j'||_{\infty}\leq G,\\ ||M_n^Tj+M_{n'}^Tj'||_{\infty}\geq ||M_{(k-1)^+}^T||_{\infty}}}T^+(j,j',n,n',x)\\
  & \:\:\: +\sum_{\substack{n,n'\in\Delta_k,\\ |n-n'|\leq 1+\log_q(1+G)}}\sum_{\substack{1\leq||j||_{\infty},||j'||_{\infty}\leq G,\\ ||M_n^Tj-M_{n'}^Tj'||_{\infty}\geq ||M_{(k-1)^+}^T||_{\infty}}}T^-(j,j',n,n',x)
 \end{aligned}
\end{equation}
and
\begin{equation}
 \begin{aligned}
  \label{144}
  S_k(x) \:\:\: = & \:\:\: \sum_{\substack{n,n'\in\Delta_k,\\ |n-n'|\leq 1+\log_q(1+G)}}\sum_{\substack{1\leq||j||_{\infty},||j'||_{\infty}\leq G,\\ 0<||M_n^Tj+M_{n'}^Tj'||_{\infty}<||M_{(k-1)^+}^T||_{\infty}}}T^+(j,j',n,n',x)\\
  & \:\:\: +\sum_{\substack{n,n'\in\Delta_k,\\ |n-n'|\leq 1+\log_q(1+G)}}\sum_{\substack{1\leq||j||_{\infty},||j'||_{\infty}\leq G,\\ 0<||M_n^Tj-M_{n'}^Tj'||_{\infty}<||M_{(k-1)^+}^T||_{\infty}}}T^-(j,j',n,n',x).
 \end{aligned}
\end{equation}
Thus we have $Y_k^2-\varsigma_k=Q_k+R_k+S_k$. In order to estimate $|\mathbb{E}[Q_k|\mathcal{F}_{k-1}]|$ we first observe that for $|n'-n|>1+\log_q(G+1)$ and $1\leq||j||_{\infty},||j'||_{\infty}\leq G$ we get
\begin{equation*}
 ||M_n^Tj\pm M_{n'}^Tj'||_{\infty}>||M_{\min(n,n')}^T||_{\infty}\geq ||M_{k^-}^T||_{\infty}.
\end{equation*}
Therefore if $C_{d,G}$ is large enough with a similar argumentation as in the proof of (\ref{115}) we get
\begin{equation*}
 \begin{aligned}
  |\mathbb{E}[Q_k|\mathcal{F}_{k-1}] |\:\:\: \leq & \:\:\: \sum_{\substack{n,n'\in\{1,\ldots, N\},\\ |n-n'|>1+\log_q(G+1)}}\sum_{1\leq||j||_{\infty},||j'||_{\infty}\leq G}\frac{1}{\lambda(\mathcal{A}_{k-1})}\\
  & \:\:\: \cdot\left|\int_{\mathcal{A}_{k-1}}T^+(j,j',n,n',x)+T^-(j,j',n,n',x)\,dx\right|\\
  \leq & \:\:\: C(5\log(G))^{2d}|\Delta_k|^2\frac{||M_{(k-1)^+}^T||_{\infty}k^{1+2\eta}2^{C'_{d,G}}}{||M_{k^-}^T||_{\infty}}\\
  \leq & \:\:\: C\eta^{-2}(5\log(G))^{2d}C_{d,G}^{2(1+\eta)}k^{2\eta}k^{1+2\eta}2^{C'_{d,G}}k^{-2(1+2\eta)}q^{-C_{d,G}}\\
  \leq & \:\:\: Ck^{-1}
 \end{aligned}
\end{equation*}
for some constant $C>0$ depending only on $q$. Thus we have
\begin{equation}
\label{124}
\left|\left|\sum_{k=1}^{K+1}\mathbb{E}[Q_k|\mathcal{F}_{k-1}]\right|\right|_2\leq C\log(K).
\end{equation}
By Lemma \ref{135} for $1\leq n,n'\leq N$ and $||j||_{\infty}\leq G$ there exists at most one $||j'||_{\infty}\leq G$ such that $||M^T_nj\pm M^T_{n'}j'||_{\infty}\leq ||M^T_{(k-1)^+}||_{\infty}$. Hence with Cauchy-Schwarz inequality we observe $||S_k||_{\infty},||\varsigma_k||_{\infty}\leq (1+\log_q(G+1))|\Delta_k|||p||_{2}^2\leq (1+\log_q(G+1))|\Delta_k|$. By definition we get $||\mathcal{R}_k||_{\infty}\leq (1+\log_q(G+1))|\Delta_k|$ and therefore we also have
\begin{equation*}
 ||R_k||_{\infty}\leq ||\mathcal{R}_k||_{\infty}+||S_k||_{\infty}+||\varsigma_k||_{\infty}\leq C(1+\log_q(G+1))|\Delta_k|.
\end{equation*}
Now we estimate $\sum_{k=1}^{K+1}\mathbb{E}[R_k|\mathcal{F}_{k-1}]$ obtaining
\begin{equation*}
 \mathbb{E}\left[\left(\sum_{k=1}^{K+1}\mathbb{E}[R_k|\mathcal{F}_{k-1}]\right)^2\right]\leq 2\mathbb{E}\left[\sum_{k,k'=1}^{K+1}\mathbb{E}[R_k|\mathcal{F}_{k-1}]\mathbb{E}[R_{k'}|\mathcal{F}_{k'-1}]\right].
\end{equation*}
For $k=k'$ we have
\begin{equation}
\label{8}
 \sum_{k=1}^{K+1}\mathbb{E}^2[R_k|\mathcal{F}_{k-1}]\leq \sum_{k=1}^{K+1}||R_k||_{\infty}^2\leq C\eta^{-2} C_{d,G}^{2(1+\eta)}(\log(G))^2K^{1+2\eta}.
\end{equation}
We may assume $k'>k$ now. Since $\mathbb{E}[R_k|\mathcal{F}_{k-1}]$ is $\mathcal{F}_{k-1}$-measurable we get
\begin{multline*}
 \left|\mathbb{E}\left[\left.\sum_{1\leq k<k'\leq K+1} \mathbb{E}[R_k|\mathcal{F}_{k-1}]\mathbb{E}[R_{k'}|\mathcal{F}_{k'-1}]\right|\mathcal{F}_{k-1}\right]\right|\\
 \begin{aligned}
 \leq & \sum_{1\leq k<k'\leq K+1}||R_k||_{\infty}|\mathbb{E}[R_{k'}|\mathcal{F}_{k-1}]|\\
 \leq & \sum_{1\leq k<k'\leq K+1}C\eta^{-1} C_{d,G}^{1+\eta}k^{\eta}(1+\log_q(G+1))|\mathbb{E}[R_{k'}|\mathcal{F}_{k-1}]|.
 \end{aligned}
\end{multline*}
Furthermore $R_{k'}$ can be represented by
\begin{equation*}
 R_{k'}(x)=\sum_{\substack{j\in\mathbb{Z}^d,\\ ||M^T_{(k'-1)^+}||_{\infty}\leq ||j||_{\infty}\leq 2G||M^T_{k'^+}||_{\infty}}}\gamma_j\cos(2\pi\langle j,x\rangle)+\delta_j\sin(2\pi\langle j,x\rangle)
\end{equation*}
where because of constraints on $n,n',j,j'$ we have 
\begin{equation*}
 \sum_{j}|\gamma_j|+|\delta_j|\leq C(5\log(G))^{2d}\eta^{-1}C_{d,G}^{1+\eta}k^{\eta}\log(G)
\end{equation*}
for some constant $C>0$ depending only on $q$.
Thus we have by using a similar argumentation as above
\begin{equation}
\begin{aligned}
\label{9}
 |\mathbb{E}[R_{k'}|\mathcal{F}_{k-1}]|\\
  \leq & \:\:\: \frac{1}{\lambda(\mathcal{A}_{k-1})}\left|\int_{\mathcal{A}_{k-1}}\sum_{j}\gamma_j\cos(2\pi\langle j,x\rangle)+\delta_j\sin(2\pi\langle j,x\rangle)\,dx\right|\\
  \leq & \:\:\: C(5\log(G))^{2d}\eta^{-1}C_{d,G}^{1+\eta}k^{\eta}\log(G)\frac{||M^T_{(k-1)^+}||_{\infty}k^{1+2\eta}2^{C'_{d,G}}}{||M^T_{(k'-1)^+}||_{\infty}}\\
 \leq & \:\:\: C(5\log(G))^{2d}\eta^{-1}C_{d,G}^{1+\eta}k^{\eta}\log(G)k^{1+2\eta}2^{C'_{d,G}}q^{-C\eta^{-1}C_{d,G}^{1+\eta}(k'-1)^{\eta}}\\
 \leq & \:\:\: Ck^{1+3\eta}q^{-(k'-1)^{\eta}}
\end{aligned}
\end{equation}
for sufficiently large $C_{d,G}$.
With (\ref{8}) and (\ref{9}) we conclude
\begin{multline}
\label{10}
 \left|\left|\sum_{k=1}^{K+1}\mathbb{E}[R_k|\mathcal{F}_{k-1}]\right|\right|_2\\
 \begin{aligned}
 \leq & \left(C\eta^{-2}C_{d,G}^{2(1+\eta)}\log(G)^2K^{1+2\eta}+\sum_{1\leq k<k'\leq K+1}k^{1+3\eta}q^{-(k'-1)^{\eta}}\right)^{1/2}\\
 \leq & C\eta^{-1}C_{d,G}^{1+\eta}\log(G)K^{1/2+\eta}.
 \end{aligned}
\end{multline}
Finally we estimate $S_k$ which can be written as
\begin{equation*}
 S_k=\sum_{0<||j||_{\infty}< ||M^T_{(k-1)^+}||_{\infty}}\gamma_j\cos(2\pi\langle j,x\rangle)+\delta_j\sin(2\pi\langle j,x\rangle),
\end{equation*}
where $\sum_{j}|\gamma_j|+|\delta_j|\leq C\eta^{-1}(5\log(G))^{2d}C_{d,G}^{1+\eta}k^{\eta}\log(G)$.
The fluctuation of $S_k$ on any atom of $\mathcal{F}_{k-1}$ is at most
\begin{equation*}
 \sum_{0<||j||_{\infty}\leq ||M^T_{(k-1)^+}||_{\infty}}(|\gamma_j|+|\delta_j|)2\pi G||M^T_{(k-1)^+}||_{\infty}\cdot 2^{-m((k-1)^+)}d\leq C\eta^{-1}k^{-(1+\eta)}.
\end{equation*}
where the inequality follows by definition of $m((k-1)^+)$ and (\ref{111}).
Therefore we get
\begin{equation*}
 |\mathbb{E}[S_k|\mathcal{F}_{k-1}]-S_k|\leq C\eta^{-1}k^{-(1+\eta)}.
\end{equation*}
Thus we have
\begin{equation}
\label{6}
 \left|\left|\sum_{k=1}^{K+1}\mathbb{E}[S_k|\mathcal{F}_{k-1}]\right|\right|_2\leq\left|\left|\sum_{k=1}^{K+1}S_k\right|\right|_2+C\eta^{-1}
\end{equation}
for some constant depending $C>0$ only on $q$.
We write
\begin{equation*}
\sum_{k=1}^{K+1}S_k(x)=\sum_{0<||j||_{\infty}< ||M^T_{K^+}||_{\infty}}\gamma'_j\cos(2\pi\langle j,x\rangle)+\delta'_j\sin(2\pi\langle j,x\rangle),
\end{equation*}
where by $L(N,G)\leq CN^{\beta}$ we have
\begin{equation}
 \label{137}
 |\gamma'_j|,|\delta'_j|\leq C\left(\eta^{-1}C_{d,G}^{1+\eta}K^{1+\eta}\right)^{\beta}.
\end{equation}
We obtain
\begin{eqnarray*}
 \left|\left|\sum_{k=1}^{K+1}S_k\right|\right|_2^2 & = & \sum_j\gamma_j'^2+\delta_j'^2\\
 & \leq & C\left(\eta^{-1}C_{d,G}^{1+\eta}K^{1+\eta}\right)^{\beta}\sum_j|\gamma_j'|+|\delta_j'|\\
 & \leq & C\left(\eta^{-1}C_{d,G}^{1+\eta}K^{1+\eta}\right)^{1+\beta}.
\end{eqnarray*}
By (\ref{6}) we have
\begin{equation}
\label{11}
 \left|\left|\sum_{k=1}^{K+1}\mathbb{E}[S_k|\mathcal{F}_{k-1}]\right|\right|_2\leq C\left(\eta^{-1}C_{d,G}^{1+\eta}K^{1+\eta}\right)^{(1+\beta)/2}+C\eta^{-1}.
\end{equation}
Using (\ref{126}), (\ref{124}), (\ref{10}) and (\ref{11}) we finally observe
\begin{equation}
\label{12}
 \mathbb{E}\left[(V_{K+1}-s_{K+1}^2)^2\right]\leq C\eta^{-2}C_{d,G}^{2(1+\eta)}\log(G)^2K^{1+\eta+\max(\eta,\beta(1+\eta))}.
\end{equation}
With Theorem \ref{5} and (\ref{38}), (\ref{13}) and (\ref{12}) we obtain
\begin{multline*}
 \sup_t\left|\mathbb{P}\left(\sum_{k=1}^{K+1}X_k\leq s_kt\right)-\Phi(t)\right|\\
 \begin{aligned}
 \leq & A\left(\frac{C\eta^{-2}C_{d,G}^{2(1+\eta)}\log(G)^2K^{1+\eta+\max(\eta,\beta(1+\eta))}}{N^2}\right)^{1/5}\\
 \leq & A\left(\frac{C\eta^{-2}C_{d,G}^{2(1+\eta)}\log(G)^2\left(C\eta^{-1/(1+\eta)}C_{d,G}^{-1}N^{1/(1+\eta)}\right)^{1+\eta+\max(\eta,\beta(1+\eta))}}{N^2}\right)^{1/5}\\
 \leq & C\eta^{-4/5}\log(G)^{2/5}C_{d,G}^{\min(1,(1-\beta)(1+\eta))/5}N^{\max(-1/(1+\eta),\beta-1)/5}.
 \end{aligned}
\end{multline*}
Now we set $\alpha=3/4$, $\eta=3/5$ and $\varepsilon=N^{-1/8}$. Thus together with (\ref{122}) we have
\begin{multline*}
 \left|\mathbb{P}\left(\sum_{n=1}^Nf(M_nx)\leq t\sigma_N\right)-\Phi(t)\right|\\
 \leq C\left(\log(d)\log(N)+\log(G)^{2/5}C_{d,G}^{1/5}\right)N^{-\min(1/8,(1-\beta)/5)}
\end{multline*}
for some $C>0$ which depends only on $q$. 
With (\ref{110}) and $N\geq Cd$ for some constant which only depends on $q$ and the constant used in (\ref{21}) we obtain
\begin{multline*}
 \left|\mathbb{P}\left(\sum_{k=1}^Nf(M_nx)\leq t\sigma_N\right)-\Phi(t)\right|\\
 \leq C\left(d^{1/5}\log(N)^{3/5}+\log(d)\log(N)\right)N^{-\min(1/8,(1-\beta)/5)}.
\end{multline*}
Therefore (\ref{139}) is proved. 

We now show (\ref{138}). Therefore we take some arbitrary large $G\in\mathbb{N}$ and repeat the proof of (\ref{139}). Observe that because of $L(N,G)=o(N)$ instead of (\ref{137}) we get
\begin{equation*}
 |\gamma'_j|,|\delta'_j|= o\left(\eta^{-1}C_{d,G}^{1+\eta}K^{2(1+\eta)}\right)
\end{equation*}
and instead of (\ref{11}) we also have
\begin{equation*}
 \left|\left|\sum_{k=1}^{K+1}\mathbb{E}[S_k|\mathcal{F}_{k-1}]\right|\right|_2= o\left(\eta^{-1}C_{d,G}^{1+\eta}K^{1+\eta}\right)+C\eta^{-1}.
\end{equation*}
With Lemma \ref{5} and (\ref{38}), (\ref{13}), (\ref{126}), (\ref{124}), (\ref{10}), (\ref{11}) and (\ref{12}) we observe
\begin{equation*}
 \sup_t\left|\mathbb{P}\left(\sum_{k=1}^{K+1}X_k\leq s_kt\right)-\Phi(t)\right|\leq h(K)
\end{equation*}
for some positive function $h$ with $\lim_{K\to\infty}h(K)=0$.
Take $\eta=3/5$ and $\varepsilon=G^{-1/6}$ Then together with (\ref{122}) for sufficiently large $N$ we observe
\begin{equation}
 \label{145}
 \left|\mathbb{P}\left(\sum_{k=1}^Nf(M_nx)\leq t\sigma_N\right)-\Phi(t)\right|\leq h(N)+CG^{-1/6}
\end{equation}
for some constant $C>0$ which depends only on $q$ and $d$. Since $G$ can be chosen arbitrary, we have shown (\ref{138}) which concludes the proof of Theorem \ref{3}.

\vspace{5ex}

\section{Proof of Theorem (\ref{16})}
\label{S32}

The following Theorem due to Strassen plays an important part in the proof of Theorem \ref{16}.

\begin{theorem}[\extcite{A10}{Lemma 2.1},\extcite{S67}{Corollary 4.5}]
\label{56}
 Let $(X_k,\mathcal{F}_k,k\geq 1)$ be a martingale differences sequence with finite fourth moments, set $V_K=\sum_{k=1}^K\mathbb{E}[X_k^2|\mathcal{F}_k]$ and assume $V_1>0$ and $V_K\to\infty$ for $K\to\infty$.
Furthermore assume
\begin{equation*}
 \lim_{K\to\infty}\frac{V_K}{b_K}=1 \quad \textnormal{a.s.}
\end{equation*}
for some sequence $(b_K)_{K\geq 1}$ of positive numbers, and
\begin{equation*}
 \sum_{K=1}^{\infty}\frac{\log(b_K)^{10}}{b_K^2}\mathbb{E}[X_K^4]<\infty.
\end{equation*}
Then we have
\begin{equation*}
 \limsup_{K\to\infty}\frac{\sum_{k=1}^KX_k}{\sqrt{2b_K\log(\log(b_K))}}=1 \quad \textnormal{a.s.}
\end{equation*}
\end{theorem}

We shall prove Theorem \ref{16} by using this result which ensures the Law of the Iterated Logarithm for a martingale differences sequence under certain conditions. Therefore we define the martingale differences sequence in the same way as in the proof of Theorem \ref{3}, i.e. $X_k=\sum_{n\in\Delta_k}\varphi_n(x)$ where the sums are taken over a certain long blocks $\Delta_k$ with small gaps between two consecutive blocks. Furthermore the functions $\varphi_n$ are piecewise constant functions which are used to approximate the trigonometric polynomials induced by the low-frequency part of $f$. Thus we need to give bounds for the remaining parts, i.e. the small blocks between two consecutive long blocks as well as the high-frequency part of $f$. Upper bounds shall be given by the following Lemma which proof is mainly based on  \cite{T62} and \cite{P75} where a similar result was obtained for the one-dimensional case.

\begin{lemma}
\label{30}
 Let $(M_n)_{n\geq 1}$ be a lacunary sequence of non-singular matrices satisfying the Hadamard gap condition (\ref{17}) . Let $f$ be a bounded periodic function of finite total variation in the sense of Hardy and Krause satisfying $\mathbb{E}[f]=0$ and $0<||f||_2\leq 1$.
Then we have
\begin{equation*}
 \limsup_{N\to\infty}\frac{\left|\sum_{n=1}^Nf(M_nx)\right|}{\sqrt{2N\log(\log(N))}}\leq C||f||^{1/4}_2 \quad\textnormal{a.e.}
\end{equation*}
for an absolute constant $C>0$ depending only on $q$.
\end{lemma}

\textit{Proof.}
For some integers $R,S\in\mathbb{N}$ set $F(R,S,x)=\left|\sum_{n=1+S}^{R+S}f(M_nx)\right|$. Furthermore for $m=\max\{l\in\mathbb{N}:2^l\leq N\}$ we have
\begin{multline}
\label{28}
\frac{F(0,N,x)}{\sqrt{2N||f||_2^{1/2}\log(\log(N))}}\\
 \begin{aligned}
  \leq & \frac{F(0,2^m,x)}{\sqrt{2\cdot2^m||f||_2^{1/2}\log(\log(2^m))}}+\sum_{l=\lceil m/3\rceil}^{m-1}\frac{F\left(2^m+\mu_{l+1}2^{l+1},\delta_l2^l,x\right)}{\sqrt{2\cdot2^m||f||_2^{1/2}\log(\log(2^m))}}\\
  & +\frac{F\left(2^m+\mu_{\lceil m/3\rceil}2^{\lceil m/3\rceil},N^*,x\right)}{\sqrt{2\cdot2^m||f||_2^{1/2}\log(\log(2^m))}}
  \end{aligned}
\end{multline}
where $\mu_l\in\{0,\ldots,2^{m-l}-1\}$ and $\delta_l\in\{0,1\}$ for all $l$ and the integer $N^*$ is given by $N^*=N-2^m-\mu_{\lceil m/3\rceil}2^{\lceil m/3\rceil}$. Let $\phi(K)=\sqrt{2K\log(\log(K))}$. Now define the sets
\begin{eqnarray*}
 D(m) & = & \left\{\left|F(0,2^m,x)\right|>16C_1||f||_2^{1/4}\phi(2^m)\right\},\\
E(m,l,\mu_{l+1}) & = & \left\{\left|F(2^m+\mu_{l+1}2^{l+1},2^l,x)\right|>16C_1\cdot 2^{(l-m)/6}||f||_2^{1/4}\phi(2^m)\right\},
\end{eqnarray*}
for some absolute constant $C_1>0$ to be specified later. We are now going to show that for any $\varepsilon>0$ there exists some $m_0\in\mathbb{N}$ such that
\begin{equation}
\label{27}
 \mathbb{P}\left(\bigcup_{m\geq m_0}\left(D(m)\cup\bigcup_{l=\lceil m/3\rceil}^{m-1}\bigcup_{\mu_{l+1}=0}^{2^{m-l-1}-1}E(m,l,\mu_{l+1})\right)\right)<\varepsilon.
\end{equation}
In order to show this inequality we apply the following inequality for suitable choices of $R$, $R'$, $S$, $Z$ and $\alpha$:
\begin{multline}
 \label{78}
 \mathbb{P}\left(\left|\sum_{n=1+S}^{R+S}f(M_nx)\right|>Z||f||^{1/4}_2\sqrt{2R\log(\log(R'))}\right)\\
\leq\frac{4C_2}{Z^2||f||^{1/2}_2R^{\alpha/2}\log(\log(R'))}+e^{-Z/4C_1\cdot||f||_2^{-1/2}\log(\log(R'))},   
\end{multline}
with $C_2=C_3(\log(d)+1)$ where (\ref{15}) is satisfied with $C=C_3$. Let $p_{R^{\alpha}}$ be the $R^{\alpha}$th Fej\'er mean of $f$ and let $r_{R^{\alpha}}=f-p_{R^{\alpha}}$. We obtain
\begin{multline}
 \label{128}
 \mathbb{P}\left(\left|\sum_{n=1+S}^{R+S}f(M_Nx)\right|>Z||f||^{1/4}_2\sqrt{2R\log(\log(R'))}\right)\\
 \begin{aligned}
  \leq & \mathbb{P}\left(\left|\sum_{n=1+S}^{R+S}p_{R^{\alpha}}(M_nx)\right|>\frac{Z}{2}||f||^{1/4}_2\sqrt{2R\log(\log(R'))}\right)\\
  & +\mathbb{P}\left(\left|\sum_{n=1+S}^{R+S}r_{R^{\alpha}}(M_nx)\right|>\frac{Z}{2}||f||^{1/4}_2\sqrt{2R\log(\log(R'))}\right).
 \end{aligned}
\end{multline}
Using Lemma \ref{94} and \ref{102} and also Chebyshev's inequality the second part can be estimated by
\begin{equation}
 \begin{aligned}
  \mathbb{P}\left(\left|\sum_{n=1+S}^{R+S}r_{R^{\alpha}}(M_nx)\right|>\frac{Z}{2}||f||^{1/4}_2\sqrt{2R\log(\log(R'))}\right) & \leq \frac{2\mathbb{E}\left[\left(\sum_{n=1+S}^{R+S}r_{R^{\alpha}}(M_nx)\right)^2\right]}{Z^2/4\cdot2||f||^{1/2}_2R\log(\log(R'))}\\
  & \leq \frac{4C_2\log(d)}{Z^2||f||^{1/2}_2\cdot R^{\alpha/2}\log(\log(R'))}.
 \end{aligned}
\end{equation}
In order to estimate the first term on the right-hand side of (\ref{128}) we apply the techniques used in the proof of Lemma \ref{123}. We shall show
\begin{equation}
 \label{129}
 \mathlarger{\mathlarger{\int}}_{[0,1)^d}\exp\left(\kappa_R\sum_{n=1+S}^{R+S}p_{R^{\alpha}}(M_nx)\right)\,dx\leq e^{C_1\kappa_{R}^2||f||_2R}
\end{equation}
for large enough $R$, suitable $\kappa_{R}$ and some absolute constant $C_1>0$.
Therefore we set
\begin{equation*}
 \kappa_R=\frac{1}{4C_1}\sqrt{\frac{2\log(\log(R'))}{R}}               
\end{equation*}
and $P=\lfloor R^{1/6\alpha}\rfloor$ for some $\alpha\geq 1$ as well as $l=\lfloor R/2P\rfloor$. Without loss of generality we may assume $\kappa_RP||f||_{\infty}\leq 1$. Again we define
\begin{equation*}
 U_m(x)=\sum_{n=Pm+1+S}^{P(m+1)+S}p_{R^{\alpha}}(x).
\end{equation*}
By Cauchy-Schwarz inequality we obtain
\begin{equation*}
 \mathlarger{\mathlarger{\int}}_{[0,1)^d}\exp\left(\kappa_R\sum_{n=1+S}^{R+S}p_{R^{\alpha}}(M_nx)\right)\,dx\leq I(\kappa_R,l)^{1/2}I'(\kappa_R,l)^{1/2}
\end{equation*}
where $I(\kappa_R,l)$ and $I'(\kappa_R,l)$ are defined similarly as in the proof of Lemma \ref{123}. 
Since $|\kappa_RU_{2m}(x)|\leq \kappa_RP\leq 1$ we estimate $I(\kappa_R,l)$ by using $e^z\leq 1+z+z^2$ which holds for $|z|\leq 1$. Thus we get
\begin{equation*}
 I(\kappa_R,l)\leq \mathlarger{\mathlarger{\int}}_{[0,1)^d}\prod_{m=0}^{l-1}(1+\kappa_RU_{2m}(x)+\kappa_R^2U_{2m}(x)^2)\,dx.
\end{equation*}
We define
\begin{equation}
 \label{150}
 \begin{aligned}
  V_{2m}(x) \:\:\: = & \:\:\: \sum_{n,n'=2Pm+1+S}^{(2m+1)P+S}\sum_{\substack{1\leq ||j||_{\infty},||j'||_{\infty}\leq R^{\alpha},\\||M^T_nj+M^T_{n'}j'||_{\infty}<||M^T_{m'}||_{\infty}}}T^+(j,j',n,n',x)\\
  & \quad\quad\quad\quad\quad\:\:\: +\sum_{\substack{1\leq ||j||_{\infty},||j'||_{\infty}\leq R^{\alpha},\\||M^T_nj-M^T_{n'}j'||_{\infty}<||M^T_{m'}||_{\infty}}}T^-(j,j',n,n',x)
 \end{aligned}
\end{equation}
where $T^+$ and $T^-$ are defined similarly as in (\ref{130}) and $m'=\lfloor 2Pm+S-\alpha\log_q(R)\rfloor$. It is easy to see that because of Lemma \ref{135} for any $n$, $n'$ and $j$ there exists at most one $j'$ such that $||M^T_nj\pm M^T_{n'}j'||_{\infty}<||M^T_{m'}||_{\infty}$. Furthermore for any $n$ and $n'$ with $||j||_{\infty},||j'||_{\infty}\leq 1/2\cdot q^{|n-n'|}$ we have
\begin{equation*}
 ||M^T_nj\pm M^T_{n'}j'||_{\infty}\geq \frac{1}{2}q^{|n-n'|}||M^T_{\min(n,n')}||_{\infty}\geq ||M^T_{m'}||_{\infty}
\end{equation*}
where the second inequality follows for some $R$ large enough. We conclude that for $||M^T_nj\pm M^T_{n'}j'||_{\infty}<||M^T_{m'}||_{\infty}$ we have $\max(||j||_{\infty},||j'||_{\infty})>1/2\cdot q^{|n-n'|}$.
Therefore by Cauchy-Schwarz inequality and an argumentation similar as in the proof of Lemma \ref{94} we observe
\begin{equation*}
 |V_{2m}(x)|\leq C\sum_{n,n'=2Pm+1+S}^{(2m+1)P+S}||p_{R^{\alpha}}||_2d^{1/2}q^{-|n-n'|/2}\leq Cd^{1/2}||f||_2P
\end{equation*}
for some constant $C>0$. Define $W_{2m}(x)=U_{2m}^2(x)-V_{2m}(x)$. Then we have
\begin{multline*}
 \kappa_RU_{2m}(x)+\kappa_R^2W_{2m}(x)\\
 =\sum_{||M^T_{m'}||_{\infty}\leq ||j||_{\infty}\leq 2R^{\alpha}||M^T_{(2m+1)P}||_{\infty}}\alpha_j\cos(2\pi\langle j,x\rangle)+\beta_j\sin(2\pi\langle j,x\rangle)
\end{multline*}
with suitable $\alpha_j$,$\beta_j$ for all $j\in\mathbb{Z}^d\backslash\{0\}$. Now for $0\leq k\leq m$ let $j_{2k}$ be any frequency vector of the trigonometric polynomial $Cd^{1/2}P+\kappa_RU_{2k}(x)+\kappa_R^2W_{2k}(x)$. If $R$ is large enough such that $P>\log_q(R^{\alpha})+\log_q(3R^{\alpha})$ then with a similar argumentation as in (\ref{131})
we get $\sum_{k=0}^mj_{2k}\neq 0$ for $j_{2m}\neq 0$. We observe
\begin{eqnarray*}
 I(\kappa_R,l) & \leq & \mathlarger{\mathlarger{\int}}_{[0,1)^d}\prod_{m=0}^{l-1}\left(1+Cd^{1/2}\kappa_R^2||f||_2P+\kappa_RU_{2m}(x)+\kappa_R^2W_{2m}(x)\right)\,dx\\
 & \leq & \mathlarger{\mathlarger{\int}}_{[0,1)^d}(1+Cd^{1/2}\kappa_R^2||f||_2P)\,dx.
\end{eqnarray*}
Therefore (\ref{129}) follows immediately with some constant $C_1>1$ which depends on $q$ and $d$.
Then by Markov's inequality we observe
\begin{multline}
 \mathbb{P}\left(\left|\sum_{n=1+S}^{R+S}p_{R^{\alpha}}(M_nx)\right|>\frac{Z}{2}||f||^{1/4}_2\sqrt{2R\log(\log(R'))}\right)\\
\begin{aligned}
&\leq 2\frac{\mathbb{E}\left[\exp\left(\kappa_R\sum_{n=1+S}^{R+S}p_{R^{\alpha}}(M_nx)\right)\right]}{\exp\left(\kappa_RZ/2\cdot||f||^{1/4}_2\sqrt{2R\log(\log(R'))}\right)}\\
&\leq 2\exp\left(C_1\kappa_R^2||f||_2R-\kappa_R\frac{Z}{2}||f||^{1/4}_2\sqrt{2R\log(\log(R'))}\right)\\
&\leq 2\exp\left(-\frac{Z\log(\log(R'))}{4C_1||f||_2^{1/2}}\right)
\end{aligned}
\end{multline}
where the last line follows by
\begin{equation*}
 C_1\kappa_R^2||f||_2R\leq \kappa_R\frac{Z}{4}||f||_2^{1/4}\sqrt{2R\log(\log(R'))}
\end{equation*}
for $||f||_2\leq 1$ and $Z\geq 1$. Now take $R=R'=2^m$, $S=0$, $Z=16C_1$ and $\alpha=2$. Then we have
\begin{equation*}
 \mathbb{P}(D(m))\leq\frac{C_2}{256C_1^2||f||^{1/2}_22^m\log(\log(2^m))}+e^{-4\cdot\log(\log(2^m))}.
\end{equation*}
It is easy to see that for any $\varepsilon>0$ there is an $m_0\in\mathbb{N}$ such that
\begin{equation}
\label{25}
 \sum_{m\geq m_0}\mathbb{P}(D(m))\leq\sum_{m\geq m_0}C\cdot 2^{-m}+\sum_{m\geq m_0}(\log(2)\cdot m)^{-4}\leq\frac{\varepsilon}{2}.
\end{equation}
To estimate $\mathbb{P}(E(m,l,\mu_{l+1}))$ we take $R=2^l$, $R'=2^m$, $S=2^m+\mu_{l+1}2^{l+1}$, $\alpha=2$ and $Z=16C_12^{(m-l)/3}$. First observe
\begin{eqnarray*}
 \mathbb{P}\left(E(m,l,\mu_{l+1})\right) & = & \mathbb{P}\left(F(2^m+\mu_{l+1}2^{l+1},2^l,x)>Z||f||^{1/4}_2\sqrt{2\cdot 2^l\log(\log(2^m))}\right)\\
& \leq & \frac{C_2||f||_2^{-1/2}}{256C_1^22^{2(m-l)/3}2^l\log(\log(2^m))}+e^{-4\cdot 2^{(m-l)/3}\log(\log(2^m))}.
\end{eqnarray*}
We have
\begin{multline}
\label{26}
 \sum_{m\geq m_0}\sum_{l=\lceil m/3\rceil}^{m-1}\sum_{\mu_{l+1}=0}^{2^{m-l-1}-1}\frac{C_2||f||_{2}^{-1/2}}{256C_1^22^{2(m-l)/3}2^l\log(\log(2^m))}\\
 \leq C\sum_{m\geq m_0}2^{m/3}\sum_{l=\lceil m/3\rceil}^{m-1}2^{-4l/3}\leq C\sum_{m\geq m_0}2^{(1/3-4/9)m}\leq \frac{\varepsilon}{4}
\end{multline}
for $m_0$ sufficiently large. Furthermore we obtain
 \begin{multline}
 \sum_{m\geq m_0}\sum_{l=\lceil m/3\rceil}^{m-1}\sum_{\mu_{l+1}=0}^{2^{m-l-1}-1}\exp\left(-4\cdot 2^{(m-l)/3}\log(\log(2^m))\right)\\
  \begin{aligned}
   \label{24} 
   = & \sum_{m\geq m_0}\sum_{l=\lceil m/3\rceil}^{m-1}\exp\left(\left(\frac{(m-l-1)\log(2)}{\log(\log(2^m))}-4\cdot 2^{(m-l)/3}\right)\log(\log(2^m))\right)\\
   \leq & \sum_{m\geq m_0}\sum_{l=\lceil m/3\rceil}^{m-1}\exp\left(-2\cdot 2^{(m-l)/3}\log(\log(2^m))\right)\\
   \leq & \sum_{m\geq m_0}\sum_{\nu=1}^{\lfloor 2m/3\rfloor}\exp\left(-2\cdot 2^{\nu/3}\log(\log(2^m))\right)\\
   \leq & \sum_{m\geq m_0}\sum_{\nu=2}^{\infty}\exp(-\nu\log(\log(2^m)))\leq C\sum_{m\geq m_0}(m\log(2))^{-2}\leq \frac{C}{m_0}\\
   \leq & \frac{\varepsilon}{4}.
  \end{aligned}
 \end{multline}
Inequalities (\ref{25}), (\ref{26}) and (\ref{24}) yield (\ref{27}). With help of this inequality we know that for any $\varepsilon>0$ there is an $N_0\in\mathbb{N}$ such that by (\ref{28}) for each $N\geq N_0$ we have
\begin{equation*}
 \frac{\left|\sum_{n=1}^Nf(M_nx)\right|}{||f||_2^{1/4}\sqrt{2N\log(\log(N))}}\leq 8C_1\left(1+\sum_{l=1}^{\infty}2^{-l/6}\right)+\frac{||f||_{\infty}}{N^{1/6}}
\end{equation*}
on a set of measure which is bounded from below by $1-\varepsilon$. Thus we obtain
\begin{equation*}
 \limsup_{N\to\infty}\frac{\left|\sum_{n=1}^Nf(M_nx)\right|}{\sqrt{2N\log(\log(N))}}\leq C||f||^{1/4}_2 \quad\textnormal{a.e.}
\end{equation*}
for some absolute constant $C$ which concludes the proof.

\vspace{5ex}

\textit{Proof of Theorem \ref{16}.} In order to proof the Theorem we repeat the prove of (\ref{12}). But here we take some arbitrary fixed integer $G$ which is sufficiently large and without loss of generality we may assume $N\geq d^{-1}G$. Observe that for fixed $G$ the definition of the blocks $\Delta'_k$ and $\Delta_k$ and therefore also the definition of the random variables $(X_k)_{k\geq 1}$ does not depend on $N$. We use $L(N,G)=\mathcal{O}(N/(\log(N))^{1+\varepsilon})$ where the implied constant may depend on $G$. Therefore instead of (\ref{12}) we get
\begin{equation*}
 ||V_{K}-s_{K}^2||_2\leq c\frac{K^{1+\eta}}{(\log(K))^{(1+\varepsilon)/2}},
\end{equation*}
where $c>0$ as in the remainder of this proof denotes a constant depending on $q$, $d$ and $G$ which may vary from line to line.
In the remainder of this proof we follow the ideas used in \cite{A10} and \cite{AFF13}.
Now we define a new probability space by taking the product of $[0,1)^d$ on which $X_k$ is defined and another probability space on which independent random variables $\xi_1,\xi_2,\ldots$ with $\mathbb{P}(\xi_n=-1)=\mathbb{P}(\xi_n=1)=1/2$ for all $n\in\mathbb{N}$ are defined. For any $k\in\mathbb{N}$ we put $\Xi_k=\sum_{n\in\Delta_k}\xi_k$.
For $m\in\mathbb{N}$ we define a martingale differences sequence $(\tilde{X}_{m,k},\tilde{\mathcal{F}}_k)$ by taking the $\sigma$-field
$\tilde{\mathcal{F}}_k=\mathcal{F}_k\times\sigma(\Xi_1,\ldots,\Xi_k)$ and setting $\tilde{X}_{m,k}=X_k+\Xi_k/m$. We further put $\tilde{s}_K^2=s_K^2+\sum_{k=1}^K|\Delta_k|/m^2$.
With (\ref{113}) we get $||X_k^4-Y_k^4||_{\infty}\leq ||X_k^2-Y_k^2||_{\infty}\cdot||X_k^2+Y_k^2||_{\infty}\leq c|\Delta_k|^4k^{-(1+2\eta)}\leq c|\Delta_k|^2$ for some constant $c>0$. By Lemma \ref{123} we have $\mathbb{E}[\tilde{X}_{m,k}^4]\leq \mathbb{E}[Y_k^4]+||Y_k^4-X_k^4||_{\infty}+ck^{\eta}\leq c|\Delta_k|^2$. Furthermore we obtain $\mathbb{E}[\tilde{X}_{m,k}^2|\tilde{\mathcal{F}}_k]=\mathbb{E}[X_{k}^2|\mathcal{F}_k]+|\Delta_k|/m^2$.
Thus we have
\begin{equation*}
 \tilde{V}_{K}=\sum_{k=1}^{K}\mathbb{E}[\tilde{X}_{m,k}^2|\tilde{\mathcal{F}}_k]=V_{K}+\frac{1}{m^2}\sum_{k=1}^{K}|\Delta_k|\geq\frac{1}{m^2}\sum_{k=1}^{K}|\Delta_k|
\end{equation*}
and
\begin{equation}
\label{44}
 ||\tilde{V}_{K}-\tilde{s}_{K}^2||_2\leq c\frac{K^{1+\eta}}{(\log(K))^{(1+\varepsilon)/2}}.
\end{equation}
We now are going to show
\begin{equation}
\label{35}
 \tilde{V}_{K}=\tilde{s}_{K}^2+o\left(\tilde{s}_{K}^2\left(\log\left(\tilde{s}_{K}^2\right)\right)^{-\varepsilon/4}\right) \quad \textnormal{a.s.}
\end{equation}
Since $\varsigma_{K}=s_{K}^2-s_{K-1}^2\leq cK^{\eta}$ we have $c_1K^{1+\eta}\leq \tilde{s}_{K}^2\leq c_2K^{1+\eta}$ for constants $c_1,c_2>0$ which depend on $q$, $d$ and $G$. We also get $s_{K'}^2-s_K^2\leq c_2(K')^{\eta}(K'-K)$ and furthermore $\tilde{s}_{K'}^2-\tilde{s}_K^2\leq c_2(K')^{\eta}(K'-K)$ for $K'\geq K$.
Set $\alpha=1-\varepsilon/2+\varepsilon^2/4$ and define $K_l=\lfloor 2^{l^\alpha}\rfloor$. Then we have
\begin{equation*}
 \frac{K_{l+1}}{K_l}=1+\mathcal{O}(l^{\alpha-1})=1+\mathcal{O}((\log(K_l))^{(\alpha-1)/\alpha})=1+o((\log(K_l))^{-\varepsilon/4})
\end{equation*}
resp.
\begin{equation*}
 |K_{l+1}-K_l|=o(K_l(\log(K_l))^{-\varepsilon/4}).
\end{equation*}
We obtain
\begin{eqnarray*}
 0 & \leq & \tilde{s}_{K_{l+1}}^2-\tilde{s}_{K_l}^2\leq c_2K_{l+1}^{\eta}(K_{l+1}-K_l)=K_{l+1}^{\eta}\cdot o\left(K_l(\log(K_l))^{-\varepsilon/4}\right)\\
& = & o\left(K_l^{1+\eta}(\log(K_l))^{-\varepsilon/4}\right)=o\left(\tilde{s}_{K_l}^2(\log(K_l))^{-\varepsilon/4}\right)
\end{eqnarray*}
or
\begin{equation*}
 \frac{\tilde{s}_{K_{l+1}}^2}{\tilde{s}_{K_l}^2}=1+o\left((\log(K_l))^{-\varepsilon/4}\right).
\end{equation*}
Since
\begin{multline*}
\sum_{l=1}^{\infty}\mathbb{P}\left(\left|\tilde{V}_{K_l}-\tilde{s}_{K_l}^2\right|>\tilde{s}_{K_l}^2(\log(K_l))^{-\varepsilon/4}\right)\\
\leq 2\sum_{l=1}^{\infty}\mathbb{E}\left[\left(\frac{\tilde{V}_{K_l}-\tilde{s}_{K_l}^2}{\tilde{s}_{K_l}^2(\log(K_l))^{-\varepsilon/4}}\right)^2\right]\leq c\sum_{l=1}^{\infty}\left(\log(K_l)\right)^{-1-\varepsilon/2}<\infty
\end{multline*}
by Borel-Cantelli-Lemma we have
\begin{equation*}
 \left|\tilde{V}_{K_l}-\tilde{s}_{K_l}^2\right|=o\left(\tilde{s}_{K_l}^2(\log(K_l))^{-\varepsilon/4}\right) \quad \textnormal{a.s.}
\end{equation*}
For $K_l\leq K<K_{l+1}$ we have
\begin{equation*}
 \left(\tilde{V}_{K_l}-\tilde{s}_{K_l}^2\right)+\left(\tilde{s}_{K_l}^2-\tilde{s}_{K_{l+1}}^2\right)\leq \tilde{V}_K-\tilde{s}_K^2\leq \left(\tilde{V}_{K_{l+1}}-\tilde{s}_{K_{l+1}}^2\right)+\left(\tilde{s}_{K_{l+1}}^2-\tilde{s}_{K_{l}}^2\right).
\end{equation*}
Therefore we have proved (\ref{35}).
By Lemma \ref{123} we get $\mathbb{E}\left[\tilde{X}_{m,K}^4\right]\leq c|\Delta_K|^2\leq cK^{2\eta}$.
Thus we have
\begin{equation}
\label{36}
 \sum_{K=1}^{\infty}\frac{(\log(\tilde{s}_K^2))^{10}}{(\tilde{s}_K^2)^2}\mathbb{E}\left[\tilde{X}_{m,K}^4\right]\leq c\sum_{K=1}^{\infty}\frac{(\log(K))^{10}}{K^2}<\infty.
\end{equation}
Now we apply Theorem \ref{56}. Therefore we get
\begin{equation}
 \label{97}
 \limsup_{K\to\infty}\frac{\left|\sum_{k=1}^K\tilde{X}_{m,k}\right|}{\sqrt{2\tilde{s}_K^2\log(\log(\tilde{s}_K^2))}}=1 \quad\textnormal{a.s.}
\end{equation}
By Lemma \ref{102} for any $N\in\mathbb{N}$ we see that
\begin{multline*}
 \left|\frac{1}{\sqrt{N}}\left|\left|\sum_{n=1}^Nf(M_nx)\right|\right|_2-\frac{1}{\sqrt{N}}\left|\left|\sum_{n=1}^Np(M_nx)\right|\right|_2\right|\\
 \leq \frac{1}{\sqrt{N}}\left|\left|\sum_{n=1}^Nr(M_nx)\right|\right|_2 \leq C(\log(d)||r||_2^2+||r||_2)^{1/2}
\end{multline*}
for some absolute constant $C>0$ which depends only on $q$. Therefore we get
\begin{multline*}
 \left|\Sigma_{f,M_n}-\limsup_{N\to\infty}\frac{1}{N}\mathlarger{\mathlarger{\int}}_{[0,1)^d}\left(\sum_{n=1}^Np(M_nx)\right)^2\,dx\right|\\
 \begin{aligned}
  \leq & C\left(\Sigma_{f,M_n}^{1/2}+(\log(d)||r||_2^2+||r||_2)^{1/2}\right)\left(\log(d)||r||_2^2+||r||_2\right)^{1/2}\\
  \leq & Cd^{1/2}\log(d)G^{-1/2}
 \end{aligned}
\end{multline*}
for some absolute constant $C>0$ which only depends on $q$.
A similar estimate holds for $\liminf$.
We obtain
\begin{eqnarray*}
 \limsup_{K(N)\to\infty}\frac{\tilde{s}_K^2\log(\log(\tilde{s}_K^2))}{s_K^2\log(\log(s_K^2))} & \leq & \frac{\limsup_{K(N)\to\infty}s_{K(N)}^2/N+1/m^2}{\liminf_{K(N)\to\infty}s_{K(N)}^2/N}\\
 & \leq & \frac{\Sigma_{f,M_n}+Cd^{1/2}\log(d)G^{-1/2}+1/m^2}{\Sigma_{f,M_n}-Cd^{1/2}\log(d)G^{-1/2}}.
\end{eqnarray*}
Hence by (\ref{97}) we have
\begin{equation}
 \label{132}
 \limsup_{K\to\infty}\frac{\left|\sum_{k=1}^K\tilde{X}_{m,k}\right|}{\sqrt{2s_K^2\log(\log(s_K^2))}}\leq\sqrt{\frac{\Sigma_{f,M_n}+Cd^{1/2}\log(d)G^{-1/2}+1/m^2}{\Sigma_{f,M_n}-Cd^{1/2}\log(d)G^{-1/2}}} \quad\textnormal{a.s.}
\end{equation}
Observe that there is a similar lower bound for the term on the left-hand side. By definition of $\tilde{X}_{m,k}$ we get
\begin{equation*}
 \left|\frac{\left|\sum_{k=1}^KX_k\right|}{\sqrt{2s_K^2\log(\log(s_K^2))}}-\frac{\left|\sum_{k=1}^K\tilde{X}_{m,k}\right|}{\sqrt{2s_K^2\log(\log(s_K^2))}}\right|\leq\frac{\left|\sum_{k=1}^K\Xi_k\right|}{m\sqrt{2s_K^2\log(\log(s_K^2))}}.
\end{equation*}
Therefore simple calculation shows that
\begin{multline}
 \label{133}
 \left|\limsup_{K(N)\to\infty}\frac{\left|\sum_{k=1}^KX_k\right|}{\sqrt{2s_K^2\log(\log(s_K^2))}}-\limsup_{K(N)\to\infty}\frac{\left|\sum_{k=1}^K\tilde{X}_{m,k}\right|}{\sqrt{2s_K^2\log(\log(s_K^2))}}\right|\\
 \leq \limsup_{K(N)\to\infty}\frac{\left|\sum_{k=1}^K\Xi_k\right|}{m\sqrt{2s_K^2\log(\log(s_K^2))}}.
\end{multline}
Since $\sum_{k=1}^K\Xi_k$ is the sum of independent random variables we have
\begin{equation}
 \label{134}
 \limsup_{K(N)\to\infty}\frac{\left|\sum_{k=1}^K\Xi_k\right|}{m\sqrt{2s_K^2\log(\log(s_K^2))}}\leq \frac{\sqrt{\Sigma_{f,M_n}+Cd^{1/2}\log(d)G^{-1/2}}}{m\sqrt{\Sigma_{f,M_n}-Cd^{1/2}\log(d)G^{-1/2}}}.
\end{equation}
Because this inequality holds for any integer $m$ by (\ref{132}), (\ref{133}) and (\ref{134}) we get
\begin{equation*}
 \begin{aligned}
  \sqrt{\frac{\Sigma_{f,M_n}-Cd^{1/2}\log(d)G^{-1/2}}{\Sigma_{f,M_n}+Cd^{1/2}\log(d)G^{-1/2}}} \:\:\: & \leq \:\:\: \limsup_{K(N)\to\infty}\frac{\left|\sum_{k=1}^KX_k\right|}{\sqrt{2s_K^2\log(\log(s_K^2))}}\\
  & \leq \:\:\: \sqrt{\frac{\Sigma_{f,M_n}+Cd^{1/2}\log(d)G^{-1/2}}{\Sigma_{f,M_n}-Cd^{1/2}\log(d)G^{-1/2}}}\quad \textnormal{a.e.}
 \end{aligned}
\end{equation*}
Observe that
\begin{equation}
\label{37}
 \left|\sum_{n=1}^Np(M_nx)\right|\leq \left|\sum_{k=1}^KX_k\right|+\left|\sum_{k=1}^K(X_k-Y_k)\right|+\left|\sum_{k=1}^K Y'_k \right|+\left|\sum_{n=\tilde{N}+1}^Np(M_nx)\right|
\end{equation}
where $K$ is defined such that $N\in\Delta'_{K+1}\cup\Delta_{K+1}$ and $\hat{N}=\sum_{k=1}^K|\Delta'_k|+|\Delta_k|$.
By using $N-\hat{N}\leq cN^{\eta/(1+\eta)}$ we have $|\sum_{k=\hat{N}+1}^Np(M_nx)|\leq cN^{\eta/(1+\eta)}$. 
By (\ref{115}) we have $|\sum_{k=1}^K(X_k-Y_k)|\leq c$.
To estimate $|\sum_{k=1}^M(Y'_k)|$ we apply Lemma \ref{30} and obtain
\begin{equation}
 \label{147}
 \left|\sum_{k=1}^K Y'_k\right|\leq C||p||^{1/4}_2\sqrt{N^{1/(1+\eta)}\log(N)\log(\log(N^{1/(1+\eta)}\log(N)))}\leq cN^{1/(2+2\eta)}\log(N)^2.
\end{equation}
Plugging these estimates into (\ref{37}) gives
\begin{equation*}
 \left|\sum_{n=1}^Np(M_nx)\right|\leq \left|\sum_{k=1}^KX_k\right|+cN^{1/(2+2\eta)}\log(N)^2.
\end{equation*}
Since $s_K^2\geq cN$ we have 
\begin{equation*}
 \frac{|\sum_{n=1}^Np(M_nx)|}{\sqrt{2s_K^2\log(\log(s_K^2))}}\leq \frac{\left|\sum_{k=1}^KX_k\right|}{\sqrt{2s_K^2\log(\log(s_K^2))}}+c\frac{\log(N)^2}{N^{\eta/(2+2\eta)}}.
\end{equation*}
It follows that
\begin{eqnarray*}
 \limsup_{N\to\infty}\frac{|\sum_{n=1}^Np(M_nx)|}{\sqrt{2s_K^2\log(\log(s_K^2))}} & \leq & \limsup_{K(N)\to\infty}\frac{\left|\sum_{k=1}^KX_k\right|}{\sqrt{2s_K^2\log(\log(s_K^2))}}+\limsup_{N\to\infty}c\frac{\log(N)^2}{N^{\eta/(2+2\eta)}}\\
 & \leq & \sqrt{\frac{\Sigma_{f,M_n}+Cd^{1/2}\log(d)G^{-1/2}}{\Sigma_{f,M_n}-Cd^{1/2}\log(d)G^{-1/2}}}\quad\textnormal{a.e.}
\end{eqnarray*}
Similar arguments yield
\begin{equation*}
 \limsup_{N\to\infty}\frac{|\sum_{n=1}^Np(M_nx)|}{\sqrt{2s_K^2\log(\log(s_K^2))}}\geq \sqrt{\frac{\Sigma_{f,M_n}-Cd^{1/2}\log(d)G^{-1/2}}{\Sigma_{f,M_n}+Cd^{1/2}\log(d)G^{-1/2}}}\quad\textnormal{a.e.}
\end{equation*}
By a similar argumentation as in (\ref{125}) we have
\begin{eqnarray*}
 \Sigma_{f,M_n}-Cd^{1/2}\log(d)G^{-1/2} & \leq & \lim_{N\to\infty}\frac{s_K^2}{N}=\lim_{N\to\infty}\frac{\sigma_N^2}{N}\\
 & \leq & \Sigma_{f,M_n}+Cd^{1/2}\log(d)G^{-1/2}
\end{eqnarray*}
and we observe
\begin{equation}
 \begin{aligned}
  \frac{\Sigma_{f,M_n}-Cd^{1/2}\log(d)G^{-1/2}}{\sqrt{\Sigma_{f,M_n}+Cd^{1/2}\log(d)G^{-1/2}}} \:\:\: & \leq \:\:\: \limsup_{N\to\infty}\frac{|\sum_{n=1}^Np(M_nx)|}{\sqrt{2N\log(\log(N))}}\\
  & \leq \:\:\: \frac{\Sigma_{f,M_n}+Cd^{1/2}\log(d)G^{-1/2}}{\sqrt{\Sigma_{f,M_n}-Cd^{1/2}\log(d)G^{-1/2}}} \quad \textnormal{a.e.}
 \end{aligned}
\end{equation}
By Lemma \ref{30} for almost any $x$ we get
\begin{equation*}
 \begin{aligned}
  \limsup_{N\to\infty}\frac{\left|\sum_{n=1}^Nf(M_nx)\right|}{\sqrt{2N\log(\log(N))}} \:\:\: & \leq \:\:\: \limsup_{N\to\infty}\frac{\left|\sum_{k=1}^Np(M_nx)\right|}{\sqrt{2N\log(\log(N))}}+\limsup_{N\to\infty}\frac{\left|\sum_{n=1}^Nr(M_nx)\right|}{\sqrt{2N\log(\log(N))}}\\
  & \leq \:\:\: \frac{\Sigma_{f,M_n}+Cd^{1/2}\log(d)G^{-1/2}}{\sqrt{\Sigma_{f,M_n}-Cd^{1/2}\log(d)G^{-1/2}}}+C||r||^{1/4}_2\\
  & \leq \:\:\: \frac{\Sigma_{f,M_n}+Cd^{1/2}\log(d)G^{-1/2}}{\sqrt{\Sigma_{f,M_n}-Cd^{1/2}\log(d)G^{-1/2}}}+C(dG^{-1})^{1/8}
 \end{aligned}
\end{equation*}
and also
\begin{equation*}
 \limsup_{N\to\infty}\frac{\left|\sum_{n=1}^Nf(M_nx)\right|}{\sqrt{2N\log(\log(N))}} \geq \frac{\Sigma_{f,M_n}-Cd^{1/2}\log(d)G^{-1/2}}{\sqrt{\Sigma_{f,M_n}+Cd^{1/2}\log(d)G^{-1/2}}}-C(dG^{-1})^{1/8}.
\end{equation*}
Since $G$ can be chosen arbitrary large, we finally observe that
\begin{equation}
 \limsup_{N\to\infty}\frac{\left|\sum_{n=1}^Nf(M_nx)\right|}{\sqrt{2N\log(\log(N))}}=\sqrt{\Sigma_{f,M_n}} \quad \textnormal{a.e.}
\end{equation}
which concludes the proof.

\vspace{5ex}

\section{Proof of Theorem (\ref{63})}
\label{S33}

The proof of this Theorem is mainly based on \cite{P75}, \cite {F08} and \cite{A10}.
We only show the Law of the Iterated Logarithm for $D^*_N$, the proof of the Law of the Iterated Logarithm for $D_N$ is essentially the same.
For some integer $h>0$ and $\beta\in [0,1)^d$ set $\beta_h$ such that $\beta_{h,i}\leq \beta_{i}<\beta_{h,i}+2^{-h}$ and $\beta_h\in B_h=\{\beta\in[0,1)^d:2^h\beta_{i}\in\{0,\ldots,2^h-1\},i=1,\ldots,d\}$.
Furthermore set $\overline{[\alpha,\beta)}=[0,\beta)\backslash [0,\alpha)$ for $\alpha,\beta\in [0,1)^d$ with $\alpha_{i}\leq \beta_{i}$ for all $i\in\{1,\ldots,d\}$. Let $f_{\beta}(x)=\mathbf{1}_{[0,\beta)}(x)-\lambda([0,\beta))$ denote the centered indicator function on $[0,\beta)$ and
$f_{\alpha,\beta}(x)=\mathbf{1}_{\overline{[\alpha,\beta)}}(x)-\lambda(\overline{[\alpha,\beta)})$ the centered indicator function on $\overline{[\alpha,\beta)}$. Now choose some arbitrary fixed integer $L>0$.
We have
\begin{equation}
 \label{67}
\begin{aligned}
D^*_N(M_1x,\ldots,M_Nx) \:\:\: & = \:\:\: \sup_{\beta\in [0,1)^d}\left|\frac{\sum_{n=1}^Nf_{\beta}(M_nx)}{N}\right|\\
 & \leq \:\:\: \max_{{\beta}_L\in B_L}\left|\frac{\sum_{n=1}^Nf_{\beta_L}(M_nx)}{N}\right|+\sup_{\beta\in [0,1)^d}\left|\frac{\sum_{n=1}^Nf_{\beta_L,\beta}(M_nx)}{N}\right|.
\end{aligned}
\end{equation}
Then, since $L$ can be chosen arbitrary large, (\ref{66}) is shown if we prove
\begin{equation}
 \label{68}
\limsup_{N\to\infty}\max_{\beta_L\in B_L}\frac{\left|\sum_{n=1}^Nf_{\beta_{L}}(M_nx)\right|}{\sqrt{2N\log(\log(N))}}=\frac{1}{2} \quad \textnormal{a.e.}
\end{equation}
and
\begin{equation}
 \label{69}
\limsup_{N\to\infty}\sup_{\beta\in[0,1)^d}\frac{\left|\sum_{n=1}^Nf_{\beta_{L},\beta}(M_nx)\right|}{\sqrt{2N\log(\log(N))}}\leq C2^{-L/8} \quad \textnormal{a.e.}
\end{equation}
for some constant $C$ depending only on $d$.
For $f_{\beta_L}$ by the second part of Lemma \ref{102} we obtain $\lim_{N\to\infty}\sigma_N^2/N=||f_{\beta_L}||_2^2=\lambda([0,\beta_L))-\lambda([0,\beta_L))^2\leq 1/4$ where equality holds for $\lambda([0,\beta_L))=1/2$. By Theorem \ref{16} we have
\begin{equation}
\begin{aligned}
 \limsup_{N\to\infty}\max_{\beta_L\in B_L}\frac{\left|\sum_{n=1}^Nf_{\beta_{L}}(M_nx)\right|}{\sqrt{2N\log(\log(N))}} \:\:\: & = \:\:\: \max_{\beta_L\in B_L}\limsup_{N\to\infty}\frac{\left|\sum_{n=1}^Nf_{\beta_{L}}(M_nx)\right|}{\sqrt{2N\log(\log(N))}}\\
 & = \:\:\: \max_{\beta_L\in B_L}||f_{\beta_L}||_2=\frac{1}{2} \quad \textnormal{a.e.}
\end{aligned}
\end{equation}
Now we are going to prove (\ref{69}). For some given $N$ we set $H=\lceil m/2+\log_2(d)\rceil$ where $m=\max\{l\in\mathbb{Z}:2^l\leq N\}$. Without loss of generality we may assume $H>L$. It is easy to see that for any $x\in [0,1)^d$ we have
\begin{equation*}
 \mathbf{1}_{[0,\beta_H)}(x)\leq\mathbf{1}_{[0,\beta)}(x)\leq \mathbf{1}_{[0,\beta_{H+1})}(x)
\end{equation*}
where for convenience we set $\beta_{H+1}$ such that $\beta_{H+1,i}=\beta_{H,i}+2^{-H}$ for all $i\in\{1,\ldots,d\}^d$. Therefore we get
\begin{equation*}
 \begin{aligned}
  \mathbf{1}_{[0,\beta_H)}(x)-\lambda([0,\beta_H))-d\cdot 2^{-H} \:\:\: & \leq \:\:\: \mathbf{1}_{[0,\beta)}(x)-\lambda([0,\beta))\\
  & \leq \:\:\: \mathbf{1}_{[0,\beta_H+2^{-H})}(x)-\lambda([0,\beta_H+2^{-H}))+d\cdot 2^{-H}.
 \end{aligned}
\end{equation*}
Thus we obtain
\begin{equation}
 \label{73}
 \left|\sum_{n=1}^N\mathbf{1}_{[0,\beta)}-\lambda([0,\beta))\right|\leq\sum_{\substack{J\subset I,\\ J \neq \emptyset}}\sum_{h\in\{L+1,\ldots,H+1\}^{|J|}}\left|\sum_{n=1}^N\varphi_{J,h}(x)\right|+d2^{-H}N.
\end{equation}
Here the sum is taken over all non-empty subsets $J$ of $I=\{1,\ldots,d\}$ and $\varphi_{J,h}$ denotes the centered indicator function on the set
\begin{equation*}
\prod_{i\in I\backslash J}[0,\beta_{L,i})\times\prod_{i\in J}[\beta_{h_{i}-1,i},\beta_{h_{i},i})
\end{equation*}
for any $h\in\{L+1,\ldots,H+1\}^{|J|}$. Now with $F_{\varphi_{J,h}}(R,S,x)=\left|\sum_{n=1+S}^{R+S}\varphi_{J,h}(M_nx)\right|$ we have
\begin{equation}
 \label{74}
 \left|\sum_{n=1}^N\varphi_{J,h}(x)\right|\leq F_{\varphi_{J,h}}(0,2^m,x)+\sum_{l=\lceil m/3\rceil}^{m-1} F_{\varphi_{J,h}}(2^m+\mu_{l+1}2^{l+1},2^l,x)+CN^{1/3}.
\end{equation}
As we shall later show the system of inequalities
\begin{equation}
 \label{75}
 \begin{aligned}
 F_{\varphi_{J,h}}(0,2^m,x) \:\:\:  & \leq \:\:\: 16C_1||\varphi_{J,h}||_2^{1/4}\sqrt{2\cdot 2^m\log(\log(2^m))},\\
 F_{\varphi_{J,h}}(2^m+\mu_{l+1}2^{l+1},2^l,x) \:\:\: & \leq \:\:\: 16C_12^{(l-m)/6}||\varphi_{J,h}||_2^{1/4}\sqrt{2\cdot 2^m\log(\log(2^m))}
 \end{aligned}
\end{equation}
holds for all $m\geq m_0$, $l\in\{\lceil m/3\rceil,\ldots,m-1\}$, $J\subset I$, $h\in\{L+1,\ldots,H+1\}^{|J|}$ with $H=\lceil m/2+\log_2(d)\rceil$ and $\beta\in[0,1)^d$ on a set of measure which is bounded from below by $1-\varepsilon$ where $\varepsilon>0$ can be chosen arbitrary and $m_0$ depends on the choice of $\varepsilon$.
Simple calculation shows
\begin{equation}
 \label{76}
 \frac{1}{2}\prod_{i\in I\backslash J}c_{i}2^{-L}\prod_{i\in J}2^{-h_{i}}\leq ||\varphi_{J,h}||_2^2\leq \prod_{i\in I\backslash J}c_{i}2^{-L}\prod_{i\in J}2^{-(h_{i}-1)}.
\end{equation}
where $c_{i}=2^L\beta_{L,i}\in\{0,\ldots,2^L-1\}$ for all $i\in I\backslash J$. Therefore we have
\begin{equation}
 \label{77}
 \begin{aligned}
  \sum_{\substack{J\subset I,\\ J\neq \emptyset}}\sum_{h\in\{L+1,\ldots,H+1\}^{|J|}}||\varphi_{J,h}||_2^{1/4} \:\:\: & \leq \:\:\: C\sum_{\substack{J\subset I,\\ J\neq \emptyset}}\sum_{h\in\{L+1,\ldots,H+1\}^{|J|}}\prod_{i\in I\backslash J}(c_{i}2^{-L})^{1/8}\prod_{i\in J}2^{-h_{i}/8}\\
  & \leq \:\:\: C\sum_{\substack{J\subset I,\\ J\neq \emptyset}}\sum_{h\in\{L+1,\ldots,H+1\}^{|J|}}\prod_{i\in J}2^{-h_{i}/8}\\
  & \leq \:\:\: C\sum_{\substack{J\subset I,\\ J\neq \emptyset}}2^{-L/8}\leq C\cdot 2^{-L/8}
 \end{aligned}
\end{equation}
where the constant $C>0$ depends only on $d$.
Furthermore we have $1+\sum_{l=1}^{\infty}2^{-l/6}\leq C$ for some absolute constant $C>0$. Combining (\ref{73}), (\ref{74}), (\ref{75}) and (\ref{77}) we finally obtain (\ref{69}). Thus it remains to show (\ref{75}).
To prove (\ref{75}) we apply the techniques used in the proof of Lemma \ref{30}. Since (\ref{75}) shall hold for any function $\varphi_{J,h}$ we first encounter all possible choices for given $J$ and $h$. We set $h'_{i}=h_{i}$ for $h_{i}\leq H$ and $h'_{i}=H$
for $h_{i}=H+1$. Therefore by definition $\varphi_{J,h}$ is a centered indicator function on a set of the form
\begin{equation*}
 \prod_{i\in I\backslash J}[0,2^{-L}c_{i})\times\prod_{i\in J}[2^{-(h_{i}-1)}a_{i},2^{-(h_{i}-1)}a_{i}+2^{-h'_{i}})
\end{equation*}
with $c\in\{1,\ldots,2^L\}^{|I\backslash J|}$ and $a\in\prod_{i\in J}\{0,\ldots,2^{h_{i}-1}-1\}$. Thus each choice of $c$ and $a$ defines a function which we denote by $\varphi_{J,h}^{(c,a)}$. 
We define the sets
\begin{equation*}
\begin{aligned}
 & D(m,J,h,c,a)\\
 & = \left\{F_{\varphi_{J,h}^{(c,a)}}(0,2^m,x)>16C_1||\varphi_{J,h}^{(c,a)}||_2^{1/4}\phi(2^m)\right\},\\
 & E(m,l,\mu_{l+1},J,h,c,a)\\
 & = \left\{F_{\varphi_{J,h}^{(c,a)}}(2^m+\mu_{l+1}2^{l+1},2^l,x)>16C_1\cdot 2^{(l-m)/6}||\varphi_{J,h}^{(c,a)}||_2^{1/4}\phi(2^m)\right\},
\end{aligned}
\end{equation*}
where $\phi(K)=\sqrt{2K\log(\log(K))}$. Now we are going to prove that for any $\varepsilon>0$ there is some integer $m_0$ such that the union of all these sets with $m\geq m_0$ has total measure which is bounded from above by $\varepsilon$.
This shall be done by estimating the measure of each set with the help of (\ref{78}). We take the same choices for $S$, $R$, $R'$ and $Z$ as in (\ref{25}), (\ref{26}) and (\ref{24}) but we take $\alpha=4d+6$. It is enough to prove
\begin{multline}
 \label{80}
\sum_{\substack{J\subset I,\\J\neq\emptyset}}\sum_{h\in\{L+1,\ldots,H+1\}^{|J|}}\sum_{c\in\{1,\ldots,2^L\}^{|I|-|J|}}\sum_{a\in\prod_{i\in J}\{0,\ldots,2^{h_{i}-1}-1\}}e^{-4\cdot 2^{(m-l)/3}||\varphi_{J,h}^{(c,a)}||_2^{-1/2}\log(\log(2^m))}\\
\leq Ce^{-4\cdot 2^{(m-l)/3}\log(\log(2^m))}
\end{multline}
for some absolute constant $C>0$ depending only on $d$ where the factor $2^{(m-l)/3}$ on the right-hand side of the first line in the case of $D(m,J,h,c,a)$ becomes $1$
and
\begin{equation}
 \label{82}
\sum_{\substack{J\subset I,\\J\neq\emptyset}}\sum_{h\in\{L+1,\ldots,H+1\}^{|J|}}\sum_{c\in\{1,\ldots,2^L\}^{|I|-|J|}}\sum_{a\in\prod_{i\in J}\{0,\ldots,2^{h_{i}-1}-1\}}||\varphi_{J,h}^{c,a}||_2^{-1/2}\leq C\cdot 2^{5/8\cdot dm}
\end{equation}
for some absolute constant $C>0$ depending only on $d$. Then the conclusion follows by similar arguments as in (\ref{25}), (\ref{26}) and (\ref{24}). Observe that with $\alpha=4d+6$ we get
\begin{equation*}
 \sum_{m\geq m_0}\sum_{l=\lceil m/3\rceil}^{m-1}\sum_{\mu_{l+1}=0}^{2^{m-l-1}-1}\frac{4C_2}{256C_1^2(2^l)^{\alpha/2}\log(\log(2^m))}\cdot C\cdot 2^{5/8\cdot dm}\leq C\sum_{m\geq m_0}2^{-dm/24}\leq\frac{\varepsilon}{4}
\end{equation*}
for $m_0$ sufficiently large and we get a similar replacement for the upper bounds on the measure of the sets $D(m,J,h,c,a)$. Without loss of generality we may assume $L$ large enough. Therefore for $h_i>L$ for all $i\in I$ by (\ref{76}) we have
\begin{multline}
 \sum_{\substack{J\subset I,\\J\neq\emptyset}}\sum_{h\in\{L+1,\ldots,H+1\}^{|J|}}\sum_{c\in\{1,\ldots,2^L\}^{|I|-|J|}}\sum_{a\in\prod_{i\in J}\{0,\ldots,2^{h_{i}-1}-1\}}e^{-4\cdot 2^{(m-l)/3}||\varphi_{J,h}^{(c,a)}||_2^{-1/2}\log(\log(2^m))}\\
\begin{aligned}
 &\leq \sum_{\substack{J\subset I,\\J\neq\emptyset}}2^{dL}\sum_{h\in\{L+1,\ldots,H+1\}^{|J|}}e^{\log(2)\cdot\sum_{i\in J}(h_{i}-1)}e^{-4\cdot 2^{(m-l)/3}\prod_{i\in J}2^{h_{i}/4}\log(\log(2^m))}\\
 &\leq \sum_{\substack{J\subset I,\\J\neq\emptyset}}2^{dL}\sum_{h\in\{L+1,\ldots,H+1\}^{|J|}}e^{-4\cdot 2^{(m-l)/3}\prod_{i\in J}2^{h_{i}/8}\log(\log(2^m))}.\\
 &\leq C\sum_{\substack{J\subset I,\\J\neq\emptyset}}2^{dL}e^{-4\cdot 2^{(m-l)/3}\prod_{i\in J}2^{L/8}\log(\log(2^m))}\\
 &\leq Ce^{-4\cdot 2^{(m-l)/3}\log(\log(2^m))}
\end{aligned}
\end{multline}
for some constant $C$ depending only on $d$ and thus (\ref{80}) is proved. Moreover we have
\begin{multline}
 \sum_{\substack{J\subset I,\\J\neq\emptyset}}\sum_{h\in\{L+1,\ldots,H+1\}^{|J|}}\sum_{c\in\{1,\ldots,2^L\}^{|I|-|J|}}\sum_{a\in\prod_{i\in J}\{0,\ldots,2^{h_{i}-1}-1\}}||\varphi_{J,h}^{c,a}||_2^{-1/2}\\
\begin{aligned}
 &\leq \sum_{\substack{J\subset I,\\J\neq\emptyset}}\sum_{h\in\{L+1,\ldots,H+1\}^{|J|}}2^{5/4\cdot (|I|-|J|)L}\prod_{i\in J}2^{5/4\cdot h_{i}}\\
 &\leq C\sum_{\substack{J\subset I,\\J\neq\emptyset}}2^{5/4\cdot (|I|-|J|)L}2^{5/4\cdot |J| H}\\
 &\leq C\cdot 2^{5/8\cdot dm}
\end{aligned}
\end{multline}
for some constant $C>0$ depending only on $d$. Observe that the last line follows by $H\leq m/2+\log_2(d)+1$. Thus (\ref{82}) is proved which finally concludes the proof of Theorem \ref{63}.

\vspace{3ex}

\small{DEPT. OF MATHEMATICS, BIELEFELD UNIV., P.O.Box 100131, 33501 Bielefeld, Germany\\
\textit{E-Mail address:} \url{tloebbe@math.uni-bielefeld.de}}


\begin{thebibliography}{99999}
\bibitem{A10} Aistleitner, C.: On the law of the iterated logarithm for the discrepancy of lacunary sequences, Trans. Amer. Math. Soc., 362, 5967-5982 (2010)
\bibitem{AB10} Aistleitner, C.: Berkes, I.: On the central limit theorem for $f(n_kx)$, Probab. Theory Relat. Fields 146, 267-289 (2010)
\bibitem{A13a} Aistleitner, C.: On the law of the iterated logarithm for the discrepancy of lacunary sequences II, Trans. Amer. Math. Soc., 365, 3713-3728 (2013)
\bibitem{AFF13} Aistleitner, C., Fukuyama, K., Furuya, Y.: Optimal bound for the discrepancies of lacunary sequences, Acta Arith. 158, 229-243 (2013)
\bibitem{BP79} Berkes, I., Philipp, W.: An a.s. invariance principle for lacunary series $f(n_kx)$, Acta. Math. Acad. Hungar. 34, 141-155 (1979)
\bibitem{CB11} Conze, J.-P., Le Borgne, S.: Limit law for some modified ergodic sums, Stoch. Dyn. 11, 107-133 (2011)
\bibitem{CBR12} Conze, J.-P., Le Borgne, S., Roger, M.: Central limit theorem for stationary products of toral automorphisms, Discrete Contin. Dyn. Syst. 32, 1597-1626 (2012)
\bibitem{DT97} Drmota, M., Tichy, R.F.: Sequences, discrepancies and applications, vol. 1651 of Lecture Notes in Mathematics, Springer, Berlin, Heidelberg, New York (1997)
\bibitem{EG55} Erd\H{o}s, P., G\'al, I.S.:  On the law of iterated logarithm, Proc. Kon. Nederl. Akad. Wetensch. 58, 65-84 (1955)
\bibitem{F08} Fukuyama, K.: The law of the iterated logarithm for the discrepancies of $\{\theta^nx\}$, Acta. Math. Hungar. 118, 155-170 (2008)
\bibitem{G66} Gaposhkin, V.F.: Lacunary series and independent functions, Russian Math. Surv. 21, 3-82 (1966)
\bibitem{G70} Gaposhkin, V.F.: The central limit theorem for some weakly dependent sequences, Theory Probab. Appl. 15, 649-666 (1970)
\bibitem{HB70} Heyde, C.C., Brown, B.M.: On the departure from normality of a certain class of martingales, Ann. Math. Stat. 41, 2161-2165 (1970)
\bibitem{I51} Izumi, S.: Notes on Fourier analysis XLIV: on the law of the iterated logarithm of some sequence of functions, J. Math. (Tokyo) 1, 1-22 (1951)
\bibitem{K46} Kac, M.: On the distribution of values of sums of the type $\sum f(2^kt)$, Ann. Math. 47, 33-49 (1946)
\bibitem{K49} Kac, M.: Probability methods in some problems  of analysis and number theory, Bull. Am. Math. Soc. 55, 641-665 (1949)
\bibitem{M50} Maruyama, G.: On an asymptotic property of a gap sequence, K\^odai Math. Sem. Rep. 2, 31-32 (1950)
\bibitem{P75} Philipp, W.: Limit theorems for lacunary series and uniform distribution mod 1, Acta Arith. 26, 241-251 (1975)
\bibitem{SZ47} Salem, R., Zygmund, A.: On lacunary trigonometric series, Proc. Nat. Acad. Sci. USA 33, 333-338 (1947)
\bibitem{SZ50} Salem, R., Zygmund, A.: La loi du logarithme it\'er\'e pour les s\'eries trigo- nom\'etriques lacunaires, Bull. Sci. Math. 74, 209-224 (1950)
\bibitem{S67} Strassen, V.: Almost sure behavior of sums of independent random variables and martingales , Fifth Berkeley Symp. Math. Stat. Prob. Vol II, Part I, 315-343 (1967)
\bibitem{T61} Takahashi, S.: A gap sequence with gaps bigger than the Hadamards, Tohoku Math. J. 13, 105-111 (1961)
\bibitem{T62} Takahashi, S.: An asymptotic property of a gap sequence, Proc. Japan Acad. 38, 101-104 (1962)
\bibitem{W59} Weiss, M.: The law of the iterated logarithm for lacunary trigonometric series, Trans. Amer. Math. Soc. 91, 444-469 (1959)
\bibitem{W16} Weyl, H.: \"Uber die Gleichverteilung von Zahlen mod. Eins, Math. Ann. 77, 313-352 (1916)
\bibitem{Z68} Zaremba, S.K.: Some applications of multidimensional integration by parts, Ann. Pol. Math. 21, 85-96 (1968)
\end{thebibliography}
\end{document}